%% file: main.tex
\documentclass[11pt,letterpaper]{amsart}
\usepackage[foot]{amsaddr}

\usepackage{mathtools}
\usepackage{amsmath}
\usepackage[T1]{fontenc}
\usepackage[latin9]{inputenc}
\usepackage{geometry}
\usepackage{wrapfig}
\usepackage{multicol}
\usepackage{graphicx}
\usepackage{soul}
\usepackage{xcolor}
\usepackage{amssymb}
\usepackage{placeins}
\usepackage{bbm}
\usepackage{cite}
\usepackage{caption}
\usepackage{enumerate}
\usepackage{afterpage}
\usepackage{enumitem}
\usepackage{bmpsize}
\usepackage{hyperref}
\usepackage{tabu}
\usepackage{enumitem}
\numberwithin{equation}{section}
\usepackage{stmaryrd}
\usepackage{tikz}
\usetikzlibrary{matrix,graphs,arrows,positioning,calc,decorations.markings,shapes.symbols}
\usetikzlibrary{calc,intersections,angles,quotes,arrows.meta}

\makeatletter
\renewcommand{\email}[2][]{%
  \ifx\emails\@empty\relax\else{\g@addto@macro\emails{,\space}}\fi%
  \@ifnotempty{#1}{\g@addto@macro\emails{\textrm{(#1)}\space}}%
  \g@addto@macro\emails{#2}%
}
\makeatother

\newtheorem{theorem}{Theorem}[section]
\newtheorem{lemma}[theorem]{Lemma}
\newtheorem{proposition}[theorem]{Proposition}

{ \theoremstyle{definition}
\newtheorem{definition}[theorem]{Definition}}
{ \theoremstyle{remark}
\newtheorem{remark}[theorem]{Remark}}

\newcommand{\im}{\mathsf{i}}
\newcommand{\Real}{\mathrm{Re}\hspace{0.5mm}}
\newcommand{\Imag}{\mathrm{Im}\hspace{0.5mm}}

\newcommand{\kgeo}{K^{\mathrm{geo}}}

\newcommand{\Leb}{\operatorname{Leb}}

\newcommand{\SFb}{S}
\newcommand{\GFb}{G}
\newcommand{\sigmaq}{\sigma}
\newcommand{\fq}{f}
\newcommand{\hq}{h}
\newcommand{\zc}{z_c}
\newcommand{\pq}{p}

\title{Uniform Convergence of Pfaffian Point Process to the Airy Line Ensemble}
\date{\today}
\author{Zhengye Zhou}

\begin{document}

\begin{abstract} We consider a family of Pfaffian Schur processes whose  first coordinate marginal relates to the half--space geometric last passage percolation.    We show that the  line ensembles corresponding to the Pfaffian Schur processes with geometric weights  converge uniformly over compact sets  to the Airy
line ensemble. By detailed asymptotic analysis of the kernels, we can verify the conditions for the finite dimensional weak convergence introduced in \cite{ED24a}. By utilizing  the tightness criteria of the line ensembles established in \cite{dimitrov2024tightness}, we can further improve the finite dimensional convergence to the uniform convergence over compact sets. Moreover, using the same methodology we also show that sequences of spiked Pfaffian Schur processes converge uniformly over compact sets to the Airy wanderer line ensembles constructed in \cite{ED24a}.
\end{abstract}

\maketitle

\tableofcontents

\input{Section1.tex}

\input{Section2.tex}
\input{Section3.tex}
\input{Section4.tex}

\subsection*{Acknowledgments} 
The author sincerely thanks Evgeni Dimitrov and Chenchen Zhao for many inspiring discussions. The author is especially grateful to the anonymous reviewers for their thorough reading and insightful feedback, which substantially strengthened the exposition and presentation of the paper. The author would also like to thank  Vincent
Zhang for writing the code used to generate some of the figures in the text.

\bibliographystyle{alpha}
\bibliography{PD}

\end{document}

%% file: Section1.tex
%
%
\section{Introduction and main results}\label{Section1}
Over the past two decades, there has been significant progress in the asymptotic analysis of half-space models within the Kardar--Parisi--Zhang (KPZ) universality class. Early developments in this direction are due to Baik and Rains, who analyzed the asymptotics of the longest increasing subsequence of random involutions and the symmetrized last passage percolation (LPP) model with geometric weights~\cite{BR01a,BR01b,BR01c}. Subsequent works have investigated a variety of half-space systems, including the last passage percolation with exponential weights in a half-quadrant \cite{BBCS18},    one-dimensional polynuclear growth model in a half-space \cite{SI04}, the log-gamma polymer in a half-space  \cite{OSZ14,IMS22,BW23}. For few of these models in half-space (essentially LPP and PNG), limiting fluctuations in the
bulk can be described analysing certain explicit formulas and they follow the Airy process,
the (conjectural) universal scaling limit of models in the KPZ class; see also \cite{DM18, Pet14, CM24, OR07, FS03, BMMK07}, for other examples of proven convergence to the Airy Process. 

A natural question arising from these results is whether such half-space models also converge to the Airy line ensemble, which is a random sequence of continuous functions that arises as a scaling limit 
in random matrix theory and in a variety of models belonging to the KPZ universality class 
\cite{Spohn,Kurt03,CorHamA,Sod15,DNV23}. 
Recent progress has established this convergence for several determinantal models:  It was established that the  last passage percolation and nonintersecting random walks  converge uniformly on compact sets to the Airy line ensemble \cite{DNV23}.
Aggarwal and Huang \cite{AH23} established a strong characterization of the Airy line ensemble and used such
characterization for the KPZ and log-gamma line ensembles, building on the tightness results
of Wu\cite{Wu23}. 
In a subsequent work, Aggarwal and Huang \cite{AH23} proved convergence of lozenge tilings of general polygonal domains to the Airy line ensemble. In this paper, we study line ensembles arising naturally in the  Pfaffian Schur process and show their convergence  to the Airy line ensemble.

%
%

\subsection{Half-space geometric last passage percolation}\label{Section1.1} It is well known that the half-space last passage percolation with geometric weights along a general down-right path is distributionally equivalent to the Pfaffian Schur process \cite{BR05,BBNV18,BBCS18} (see \cite[Theorem 2.7]{DY25b}). Motivated by this correspondence, we begin by introducing the \emph{symmetrized} (\emph{half-space}) geometric last passage percolation (LPP) model.

The half-space model depends on a sequence of real parameters $\left\{a_i\right\}_{i \geq 1}$ and a parameter $c$, satisfying
\begin{equation}
a_i \geq 0, \qquad c \geq 0, \qquad a_i a_j \in [0,1), \qquad c a_i \in [0,1).    
\end{equation}
Let $W=\left(w_{i, j}: i, j \geq 1\right)$ be a random weight array defined as follows. The collection $\{w_{i,j} :i, j \geq 1\}$ consists of independent geometric random variables with
\begin{equation}\label{eq: weights}
   w_{i, j} \sim \operatorname{Geom}\left(a_i a_j\right) \text{ when } i \neq j, w_{i, i} \sim \operatorname{Geom}\left(c a_i\right),\text{  and } w_{i, j}=w_{j, i} \text{  for all } i, j \geq 1. 
\end{equation}
Here,  $X \sim \mathrm{Geom}(\alpha)$  means that $\mathbb{P}(X = k ) = \alpha^k (1- \alpha)$ for $k \in \mathbb{Z}_{\geq 0}$. We  associate the weight $w_{i,j}$ with the lattice site $(i,j) \in \mathbb{Z}^2$, see  Figure \ref{Fig.Grid}.

An {\em up-right path} $\pi$ in $\mathbb{Z}^2$ is a (possibly empty) sequence of vertices $\pi = (v_1, \dots, v_r)$ with $v_i \in \mathbb{Z}^2$, and $v_i - v_{i-1} \in \{(0,1), (1,0)\}$. For an up-right path $\pi$ in $\mathbb{Z}_{\geq 1}^2$, we define its {\em weight} by
\begin{equation}\label{Eq.PathWeight}
W(\pi) = \sum_{v \in \pi} w_v,
\end{equation}
and for any $(m,n) \in \mathbb{Z}_{\geq 1}^2$, we define the {\em last passage time} $G_1(m,n)$ by
\begin{equation}\label{Eq.LPT}
G_1(m,n) = \max_{\pi} W(\pi),
\end{equation} 
where the maximum is over all up-right paths from $(1,1)$ to $(m,n)$.

The last passage time  (\ref{Eq.LPT}) has been extensively studied in the literature. In a series of works \cite{BR01a,BR01b,BR01c},  Baik and Rains established that the fluctuations of $G_1(n,n)$ with homogeneous $a_i$, as $n \to \infty$, converge to the GOE Tracy-Widom distribution $F_{\mathrm{GOE}}$ when $c = 1$, and to the GSE Tracy-Widom distribution $F_{\mathrm{GSE}}$ when $c \in (0,1)$. Precise definitions of~$F_{\mathrm{GOE}}$ and~$F_{\mathrm{GSE}}$ can be found in~\cite{TW05}. Furthermore, Baik and Rains~\cite{BR01c} showed that when the boundary parameter is jointly scaled as~$c = 1 - \varpi \alpha_q n^{-1/3}$, the centered and scaled last passage time converges to a one-parameter family of crossover distributions~$F_{\mathrm{cross}}(\cdot; \varpi)$ interpolating between~$F_{\mathrm{GOE}}$ at~$\varpi = 0$,~$F_{\mathrm{GSE}}$ as~$\varpi \to \infty$, and the Gaussian distribution~$\Phi$ as~$\varpi \to -\infty$. 

In a related direction, Imamura and Sasamoto~\cite{SI04} studied the \emph{polynuclear growth (PNG)} model for the special choices~$c = 0$ and~$c = 1$. They proved that, for any fixed~$\kappa \in (0,1)$, the finite-dimensional distributions of~$G_1(m,n)$ in a neighborhood of~$(\kappa n, n)$, under suitable scaling, converge to those of the Airy process. In \cite{BBNV18}, Betea, Bouttier, Nejjar, and Vuleti{\'c} extended these results by computing the limiting finite-dimensional distributions of $G_1(m,n)$ near $(n,n)$ when $c$ is scaled as $c = 1 - \varpi \alpha_q n^{-1/3}$ jointly with $n$. They obtained a one-parameter family of Fredholm Pfaffians indexed by $\varpi$, which interpolate between the two cases identified in \cite{SI04}: $\varpi = 0$ corresponding to $c = 1$ and $\varpi \to \infty$ corresponding to $c = 0$. The same family of Fredholm Pfaffians was independently derived by Baik, Barraquand, Corwin, and Suidan in the context of half-space last passage percolation with exponential weights~\cite{BBCS18}, where they established a phase transition in the one-point fluctuations and identified a two-dimensional crossover among the GUE, GOE, and GSE distributions. Related crossover phenomena also arise in the facilitated exclusion process~\cite{BBCS16}, and the corresponding phase transition in the one-point fluctuations is analogous to that observed in spiked covariance matrix models~\cite{BBP05}. \\

\begin{figure}[ht]
\centering
\begin{tikzpicture}[scale=0.9]

\def\m{7} 
\def\n{5} 

\definecolor{Bg}{gray}{1.0}        
\definecolor{C2}{gray}{0.6} 

\begin{scope}[shift={(0,0)}]

  \foreach \j in {1,...,\m}{
      \node at (\j+0.5, -0.5) {\(\j\)};
  }
  \foreach \j in {1,...,\n}{
      \node at (0.35,-0.5 + \j) {\(\j\)};
  }

  \foreach \x/\y in {1/1, 1/2, 1/3, 2/3, 3/3, 3/4, 3/5, 4/5, 5/5, 6/5, 7/5} {
    \fill[C2] (\x,\y) rectangle ++(1,-1);
  }

  \foreach \i in {1,...,\n}{
    \foreach \j in {1,...,\m}{
      \draw[black] (\j,\i) rectangle ++(1,-1);
      \node at (\j+0.5, \i-0.5) {$w_{\j,\i}$};
    }
  }

  \node at (\m+1.5,-0.5) {$i$};
  \node at (0.35,-0.5 + \n + 1) {$j$};

\end{scope}

\end{tikzpicture}
\caption{The array $W = (w_{i,j}: i,j \geq 1)$ with symmetric $w_{i,j}$ as defined in \eqref{eq: weights}  and an up-right path $\pi$ (in gray) that connects $(1,1)$ to $(7,5)$.} \label{Fig.Grid}
\end{figure}

In this paper, we extend the analysis to the  higher-rank last passage times  $G(m,n) = (G_k(m,n): k \geq 1)$ defined  as follows. For $k = 1, \dots, \min(m,n)$, we define
\begin{equation}\label{Eq.HRLPT}
G_k(m,n) = \max_{\pi_1, \dots, \pi_k} \left[ W(\pi_1) + \cdots + W(\pi_k) \right],
\end{equation}
where the maximum is over $k$-tuples of pairwise disjoint up-right paths $(\pi_1, \dots, \pi_k)$ with $\pi_i$ connecting the points $(1,i)$ with $(m, n-k + i)$. When $k \geq \min(m,n) + 1$, $k$ disjoint paths as above don't exist,  so the convention is to set
\begin{equation}\label{Eq.HRLPT2}
G_k(m,n) = \sum_{i = 1}^m \sum_{j = 1}^n w_{i,j} \mbox{ for } k \geq \min(m,n) + 1.
\end{equation}

 In fact, from \cite[(2.12)]{DY25b}, we know that the sequence $\lambda(m,n) = (\lambda_k(m,n): k \geq 1)$ defined by
\begin{equation}\label{Eq.LPTLambdas}
\lambda_1(m,n) = G_1(m,n) \mbox{ and } \lambda_k(m,n) = G_k(m,n) - G_{k-1}(m,n) \mbox{ for }k\geq 2,
\end{equation}
 is a {\em partition} (i.e. a decreasing sequence of non-negative integers that is eventually zero). In addition, fix $n$ and let $m$ vary from $1$ to $n$, these partitions {\em interlace} as follows:
\begin{equation}\label{Eq.InterlaceIntro}
\lambda_1(m,n) \geq \lambda_{1}(m-1,n) \geq \lambda_{2}(m,n) \geq  \lambda_{2}(m-1,n) \geq \cdots, 
\end{equation}
where $\lambda(0,n)$ is the empty partition. By performing linear interpolation in \(m\) between consecutive points \((m, \lambda_i(m,n))\), we obtain continuous curves in \(m\), and thus view \(\{\lambda_i(m,n): i \ge 1,\, 0 \le m \le n\}\) as a sequence of random continuous functions (i.e., as  a {\em line ensemble}).

\begin{figure}[ht]
    \centering
    \begin{tikzpicture}
        \node[anchor=south west, inner sep=0] (img1) at (0,0)
            {\includegraphics[width=0.4\textwidth]{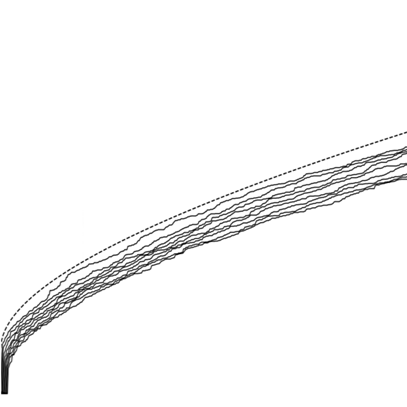}};  

        \draw[->][gray] (0.03,0.29) -- (7.2,0.29);
        \draw[->][gray] (0.03,0.29) -- (0.03,6.5);

        \draw (7.3,0.29) node[below = 0pt]{$m$};
        \fill (6.64, 0.29) circle (1.5pt) node[below=0pt] {$n$};
        \fill (0.03, 0.29) circle (1.5pt) node[below=0pt] {$0$};

        \fill (0.03, 2.37) circle (1.5pt) node[left=0pt] {$n$};
        \fill (0.03, 4.45) circle (1.5pt) node[left=0pt] {$2n$};

        \draw[black] (0.5,5.5) rectangle ++(3.5,0.6);
        \draw (2.8, 5.8) node{$nh  (m/n)$};
        \draw[-][line width = 1pt, dash pattern=on 2pt off 1pt][black] (0.6,5.8) -- (1.6,5.8);
    \end{tikzpicture}
    \caption{The top several curves of $\{\lambda_i(m,n): i \geq 1, 0 \leq m \leq n\}$ for $q = 0.5$, $c=0.8$ and $n = 500$ concentrate around the function $n h  (m/n)$ on $[0,n]$, where $h  (\kappa) = \frac{ q (q + 2 \sqrt{\kappa} + q \kappa)}{1-q^2}$.}
    \label{Fig.Simulation}
\end{figure}

The goal of this paper is to provide a detailed description of the behavior of the line ensembles constructed above. As shown by the results in the next section, under the assumptions  $a_i=q\in(0,1)$, $\kappa \in (0,1)$,   when $c \in (0,\zc)$ with $\zc$ given in (\ref{Eq.ConstBotIntro}),  the line ensembles $\{\lambda_{i}(t,n)\}_{i \geq 1}$ for $t$ near $\kappa n$ converge (after a suitable shift and scaling) to the Airy line ensemble. See Figure~\ref{Fig.Simulation} for a simulation.  Furthermore,  when $c= \zc-\sigmaq^{-1}\varpi n^{-1/3}$ with $\sigmaq$ given in (\ref{Eq.ConstBotIntro}), the line ensembles $\{\lambda_{i}(t,n)\}_{i \geq 1}$ for $t$ near $\kappa n$ instead converge (after a suitable shift and scaling) to an Airy wanderer line ensemble from \cite{AFM10,ED24a}. In the joint work with Dimitrov \cite{DZ25}, the author  derived the behavior of the line ensembles when  $c>\zc$. In this case, the top curve $\lambda_1(t,n)$ separates from the rest and behaves like a Brownian motion when appropriately shifted, and scaled horizontally by $n$ and vertically by $n^{1/2}$. On the other hand, the remaining curves $\{\lambda_{i}(t,n)\}_{i \geq 2}$ converge for $t$ near $\kappa n$ for each $\kappa \in (\kappa_0,1)$, under an appropriate shift, and $n^{2/3}$ horizontal and $n^{1/3}$ vertical scaling, to the Airy line ensemble. Furthermore, \cite{DY25} constructed the {\em half-space Airy line ensembles} as weak limits of the line ensembles described above with $\kappa = 1$ and $c = 1 - \varpi \alpha_q n^{-1/3}$. These ensembles serve as half-space analogues of the Airy line ensemble constructed in~\cite{CorHamA}.

\subsection{Main results}
In this section we present the main results of the paper. We first introduce the scaling of our line ensembles in the following definition.
\begin{definition}\label{Def.CriticalScaledLPP}
Fix $q \in (0,1)$, $\kappa \in (0, 1)$, and set
\begin{equation}\label{Eq.ConstBotIntro}
\begin{split}
&h  = \frac{ q (q + 2 \sqrt{\kappa} + q \kappa)}{1-q^2}, \hspace{2mm} p   = \frac{d}{d\kappa} h  (\kappa) = \frac{q (1 + q \sqrt{\kappa})  }{(1-q^2) \sqrt{\kappa} }, \hspace{2mm} \\
& f  = \frac{q^{1/3}}{2 \kappa^{2/3}(1-q^2)^{2/3} (q + \sqrt{\kappa})^{1/3} (1 + q \sqrt{\kappa})^{1/3}},\hspace{2mm} \zc = \frac{1 + q \sqrt{\kappa}}{q + \sqrt{\kappa}},\\
&\sigmaq = \frac{q^{1/3} (q + \sqrt{\kappa})^{5/3} }{\kappa^{1/6} (1 - q^2)^{2/3}(1+q \sqrt{\kappa})^{1/3}}.
\end{split}
\end{equation}
Fix $N \in \mathbb{N}$, set $C  _N = \lfloor h  N \rfloor$, 
and let $\lambda(m,N)$ be as in (\ref{Eq.LPTLambdas}). For $i \geq 1$, define 
\begin{equation}\label{Eq.ExtPaths}
U^{N}_i(s) = \begin{cases} \lambda_i(s,N)  &\mbox{ if } s= 1,\dots, N, \\  0 &\mbox{ if  } s \leq 0, \\ \lambda_i(N,N) &\mbox{ if } s \geq N+1. \end{cases}
\end{equation}
By linearly interpolating $(s,U^{N}_i(s))$, we may view $U_i^{N}$ as random continuous functions on $\mathbb{R}$. We define the rescaled processes $\{\mathcal{U}_i^{N} \}_{i \geq 1}$ on $\mathbb{R}$ by
$$\mathcal{U}_i^{N}(t) = [p   (1 + p  )]^{-1/2} N^{-1/3} \left( U_{i}^{N}(\lfloor \kappa N \rfloor + t N^{2/3}) - C  _N -  p   t N^{2/3}   \right) \mbox{ for } t \in \mathbb{R}.$$
We denote by $\mathcal{U}^{N} = \{\mathcal{U}^{N}_i\}_{i \geq 1}$ the corresponding line ensemble.
\end{definition}
To state the  main results of the paper we introduce the Airy wanderer line ensemble from \cite{AFM10,ED24a} in the following definition.

\begin{definition}\label{Def.AiryLE}  Given finite $J_1, J_2 \in \mathbb{Z}_{\geq 0}$ and parameters $\mathfrak{A}=\left\{a_i\right\}_{i=1}^{J_1}, \mathfrak{B}=\left\{b_i\right\}_{j=1}^{J_2}$ such that
$$
\underline{a}:=\min \left(a_i: i=1, \ldots, J_1\right)>\max \left(b_j: j=1, \ldots, J_2\right)=: \bar{b},
$$ The {\em Airy wanderer line ensemble} $\mathcal{A}^{\mathfrak{A},\mathfrak{B}} = \{\mathcal{A}_i^{\mathfrak{A},\mathfrak{B}}\}_{i \geq 1}$ is a sequence of real-valued random continuous functions, defined on $\mathbb{R}$, and are strictly ordered in the sense that $\mathcal{A}_i^{\mathfrak{A},\mathfrak{B}}(t) > \mathcal{A}_{i+1}^{\mathfrak{A},\mathfrak{B}}(t)$ for all $i \geq 1$ and $t \in \mathbb{R}$. It is uniquely specified by the following conditions. If one fixes a finite set $\mathsf{S} = \{s_1, \dots, s_m\} \subset \mathbb{R}$ with $s_1 < \cdots < s_m$, then the random measure on $\mathbb{R}^2$, defined by
\begin{equation}\label{Eq.RMS1}
M(A) = \sum_{i \geq 1} \sum_{j = 1}^m {\bf 1} \left\{\left(s_j, \mathcal{A}_{i}^{\mathfrak{A},\mathfrak{B}}(s_j) \right) \in A \right\},
\end{equation}
is a determinantal point process on $\mathbb{R}^2$ with reference measure $\mu_{\mathsf{S}} \times \mathrm{Leb}$, where $\mu_{\mathsf{S}}$ is the counting measure on $\mathsf{S}$, and $\mathrm{Leb}$ is the usual Lebesgue measure on $\mathbb{R}$, and whose correlation kernel is the {\em  the Airy wanderer kernel}   defined for $x_1, x_2 \in \mathbb{R}$ and $t_1, t_2 \in \mathsf{S}$ by 
\begin{equation}\label{Eq.S1AiryKer}
\begin{aligned}
& K_{\mathfrak{A},\mathfrak{B}}^{\mathrm{Airy}}\left( t_1, x_1 ;  t_2, x_2\right)=-\frac{\mathbf{1}\left\{ t_2> t_1\right\}}{\sqrt{4 \pi\left( t_2- t_1\right)}} \cdot e^{-\frac{\left(x_2-x_1\right)^2}{4\left( t_2- t_1\right)}-\frac{\left( t_2- t_1\right)\left(x_2+x_1\right)}{2}+\frac{\left( t_2- t_1\right)^3}{12}} \\
& +\frac{1}{(2 \pi \mathrm{i})^2} \int_{\mathcal{C}_{\gamma}^{\pi/3}} d z \int_{\mathcal{C}_{\beta}^{2\pi/3}} d w \frac{e^{z^3 / 3-x_1 z-w^3 / 3+x_2 w}}{z+ t_1-w- t_2} \cdot \prod_{i=1}^{J_1} \frac{w+ t_2-a_i}{z+ t_1-a_i} \cdot \prod_{j=1}^{J_2} \frac{z+ t_1-b_j}{w+ t_2-b_j}.
\end{aligned}
\end{equation}
In (\ref{Eq.S1AiryKer}) we have that $ \gamma+ t_1>\beta+ t_2$, in addition, $\gamma+t_1<\underline{a}$ and $\beta+t_2>\bar{b}$, and $\mathcal{C}_{z}^{\varphi}=\{z+|s|e^{\mathrm{sgn}(s)\im\varphi}, s\in \mathbb{R}\}$ is the infinite curve oriented from $z+\infty e^{-\im\varphi}$ to $z+\infty e^{\im\varphi}$. 
\end{definition}
\begin{remark}\label{Rem.AiryLE0}
When $J_1=0$, we assume $\underline{a}:=\infty$ and when $J_2=0$, we assume $\bar{b}=-\infty$.
\end{remark}

\begin{remark}\label{Rem.AiryLE1}
When $J_1=J_2=0$,  $K_{\mathfrak{A},\mathfrak{B}}^{\text {Airy }}$ reduces to the extended Airy kernel 
\begin{equation}\label{Eq.S1AiryKer2}
\begin{split}
K^{\mathrm{Airy}}(t_1,x_1; t_2,x_2) = & -  \frac{{\bf 1}\{ t_2 > t_1\} }{\sqrt{4\pi (t_2 - t_1)}} \cdot e^{ - \frac{(x_2 - x_1)^2}{4(t_2 - t_1)} - \frac{(t_2 - t_1)(x_2 + x_1)}{2} + \frac{(t_2 - t_1)^3}{12} } \\
& + \frac{1}{(2\pi \im)^2} \int_{\mathcal{C}_{\gamma}^{\pi/3}} d z \int_{\mathcal{C}_{\beta}^{2\pi/3}} dw \frac{e^{z^3/3 -x_1z - w^3/3 + x_2w}}{z + t_1 - w - t_2}
\end{split}
\end{equation}as in \cite[Definition 1.3]{DZ25}, which originates from \cite[Proposition 4.7 and (11)]{BK08}. In addition, the Airy wanderer line ensemble reduces to the Airy line ensemble $\mathcal{A} = \{\mathcal{A}_i\}_{i \geq 1}$ as in \cite[Definition 1.3]{DZ25}, whose notation we adopt throughout.
\end{remark}

\begin{remark}\label{Rem.AiryLE2}
The Airy wanderer line ensemble we defined is a special case of  the line ensemble $\left(\sqrt{2} \cdot \mathcal{L}_i^{a, b, c}(t)+t^2: i \geq 1, t \in \mathbb{R}\right)$ introduced in \cite[Theorem 1.10]{ED24a} with the parameters  $a^-=(a_1,a_2,\cdots,a_{J_1},0,\cdots)$, $b^-=(-b_1,-b_2,\cdots,-b_{J_2},0,\cdots)$ and $c^{\pm}=a_i^+=b_i^+=0$.
\end{remark}

\begin{remark}\label{Rem.AiryLE3}The extended Airy kernel is often attributed to \cite{Spohn, Kurt03}, where it arises in the study of the polynuclear growth model. A continuous version of the Airy line ensemble was subsequently constructed in \cite{CorHamA} as the weak limit of Brownian watermelons. In \cite{BP08, AFM10}, the authors introduced multi-parameter extensions of the extended Airy kernel, obtained as scaling limits of the exponential last-passage percolation model with defective rows and columns, and of the Brownian watermelon model with finitely many outliers, respectively. This family of kernels was later further generalized to the Airy wanderer kernel by Dimitrov~\cite{ED24a}, providing an infinite-parameter analogue of the extended Airy kernel.
\end{remark}

The main results of the paper, Theorem \ref{Thm.Main1} and \ref{Thm.Main2} (which are special cases of the more general Theorem
\ref{Thm.generalMain} in the main text), establish the convergence of the line ensemble $\mathcal{U}^{N}$ defined in Definition \ref{Def.CriticalScaledLPP}. Compared with previous works, which primarily establish convergence of finite-dimensional distributions via explicit formulas, our result proves \emph{uniform convergence} of the entire line ensemble, providing a stronger convergence within the same universality class.

 We denote by $C\left(\mathbb{N} \times \mathbb{R} \right)$ the space of real-valued continuous functions on $\mathbb{N} \times \mathbb{R}$ with the topology of uniform convergence over compact sets.
\begin{theorem}\label{Thm.Main1} Assume the same notation as in Definition \ref{Def.CriticalScaledLPP}, further assume $c \in (0, \zc)$ and $a_i=q$ for all $i\ge 1$. Then, $\mathcal{U}^{N}  \Rightarrow \mathcal{U}^{\infty}$ in $C\left(\mathbb{N} \times \mathbb{R} \right)$. Here, $\mathcal{U}^{\infty} = \{\mathcal{U}^{\infty}_i\}_{i \geq 1}$ is the line ensemble, defined by
$$\mathcal{U}_i^{\infty}(t) =  (2f )^{-1/2} \cdot \left( \mathcal{A}_i(f  t) - f ^2 t^2 \right) \mbox{ for } i \geq 1, t \in \mathbb{R},$$
where $\mathcal{A} = \{\mathcal{A}_i\}_{i \geq 1}$ is the Airy line ensemble from Definition \ref{Def.AiryLE}.
\end{theorem}

\begin{theorem}\label{Thm.Main2} Assume the same notation as in Definition \ref{Def.CriticalScaledLPP}, further assume $c= \zc-\sigmaq^{-1}\varpi N^{-1/3}$ with $\varpi\in \mathbb{R}$  and $a_i=q$ for all $i\ge 1$. Then, $\mathcal{U}^{N}  \Rightarrow \mathcal{U}^{\infty}$ in $C\left(\mathbb{N} \times \mathbb{R} \right)$. Here, $\mathcal{U}^{\infty} = \{\mathcal{U}^{\infty}_i\}_{i \geq 1}$ is the line ensemble, defined by
$$\mathcal{U}_i^{\infty}(t) =  (2f )^{-1/2} \cdot \left( \mathcal{A}_i^{\{\varpi\},\emptyset}(f  t) - f ^2 t^2 \right) \mbox{ for } i \geq 1, t \in \mathbb{R},$$
where $\mathcal{A}^{\{\varpi\},\emptyset} = \left\{\mathcal{A}_i^{\{\varpi\},\emptyset}\right\}_{i \geq 1}$ is the Airy wanderer line ensemble from Definition \ref{Def.AiryLE}.
\end{theorem}
\begin{remark}
We mention here that Theorems \ref{Thm.Main1} and \ref{Thm.Main2} are special cases of the more general Theorem \ref{Thm.generalMain} stated in the main text. In particular, Theorem \ref{Thm.generalMain} extends Theorems \ref{Thm.Main1} and \ref{Thm.Main2} by allowing spiked bulk parameters, leading to Airy wanderer line ensembles with additional parameters in the scaling limit.  
\end{remark}
%
%

\subsection{Main ideas and paper outline}\label{Section1.3} Our analysis starts from the distributional equality between half-space LPP and the Pfaffian Schur processes from \cite{BR05}. The definition of the Pfaffian Schur processes is recalled in Section \ref{Section2.1} and the distributional equality is detailed in Proposition \ref{Prop.LPPandSchur}. One key advantage of the Pfaffian Schur process is that it is a Pfaffian point process with an explicit correlation kernel, expressible as a double contour integral. The general formula for the kernel is presented in Proposition \ref{Prop.CorrKernel1}. In addition, we derive two alternative formulas for the correlation kernel in Lemmas \ref{Lem.PrelimitKernelsub} and \ref{Lem.PrelimitKernelscrtical}, which are suitable for asymptotic analysis in the two scaling regimes $c \in (0, \zc)$ and $c= \zc-\sigmaq^{-1}\varpi N^{-1/3}$, respectively.

In Sections \ref{Section3}  we establish the uniform convergence over compact sets of our correlation kernels from Lemmas \ref{Lem.PrelimitKernelsub} and \ref{Lem.PrelimitKernelscrtical}. The precise statements are given in Propositions \ref{Prop.KernelConvsub} and \ref{Prop.KernelConvcritical}, and proved using the method of steepest descent. One of the main technical contributions and challenges of this work lies in the detailed asymptotic analysis of the correlation kernels, in particular in identifying contour choices that remain valid across the entire range of parameters $q, c,$ and $\kappa$.

The second key feature of the Pfaffian Schur process is that it possesses the structure of a {\em geometric line ensemble} satisfying the {\em interlacing Gibbs property}; we refer the reader to \cite{dimitrov2024tightness} and \cite[Section~6]{DZ25} for an introduction and detailed discussion. This property allows the proofs of Theorem~\ref{Thm.generalMain} (and hence Theorems~\ref{Thm.Main1} and~\ref{Thm.Main2}) to be reduced to establishing finite-dimensional convergence, using the tightness framework for interlacing geometric line ensembles  developed in~\cite{dimitrov2024tightness}. The finite-dimensional convergence of the line ensembles is proved in Proposition~\ref{Prop.FinitedimBulk}, from which Theorem~\ref{Thm.generalMain} then follows directly in Section~\ref{Section4.2} by verifying the tightness criterion of~\cite{dimitrov2024tightness}.

%% file: Section2.tex
%
%
\section{Pfaffian Schur processes}\label{Section2} In this section, we introduce a family of {\em Pfaffian Schur processes} and discuss several of their key properties. The models considered here are special cases of those originally defined in~\cite{BR05}. For additional background and further developments on Pfaffian Schur processes, we refer the reader to~\cite{BR05, BBNV18, BBCS18}.

%
%
\subsection{Definitions}\label{Section2.1} Recall that a {\em partition} is a nonincreasing sequence of nonnegative integers $\lambda=\left(\lambda_1 \geqslant \lambda_2 \geqslant \cdots \geqslant \lambda_k \geqslant 0\right)$, with finitely many non-zero components. Denote by $|\lambda|=\lambda_1+\lambda_2+\cdots$ the {\em weight} of partition $\lambda$. We denote by $\emptyset$ the partition of weight $0$. Given two partitions $\lambda$ and $\mu$, we say that they {\em interlace}, denoted by $\lambda \succeq \mu$ or $\mu \preceq \lambda$, if $\lambda_1 \geq \mu_1 \geq \lambda_2 \geq \mu_2 \geq \cdots$.

Given finitely many variables  $x_1, \ldots, x_n$, the {\em skew Schur polynomials} for two partitions $\lambda, \mu$ are defined via 
\begin{equation}\label{Eq: skewshur}
  s_{\lambda / \mu}\left(x_1, \ldots, x_n\right)=\sum_{\mu=\lambda^0 \preceq \lambda^1 \preceq \cdots \preceq \lambda^n=\lambda} \prod_{i=1}^n x_i^{\left|\lambda^i-\lambda^{i-1}\right|} .  \end{equation}
For a partition $\lambda$, we also define the {\em boundary monomial} in a single variable $c$ by
\begin{equation}\label{Eq.Tau}
\tau_{\lambda}(c) = c^{\sum_{j = 1}^{\infty} (-1)^{j-1}\lambda_j} = c^{\lambda_1 - \lambda_2 + \lambda_3 - \lambda_4 + \cdots}.
\end{equation}
With the above notation in place we can define our main object of interest.
\begin{definition}\label{Def.SchurProcess} Fix $N \in \mathbb{N}$ and $c, a_1, \dots, a_N \in [0,\infty)$, such that $ca_i, a_j a_k < 1$ for $i,j,k \in \{1, \dots, N\}$ with $j \neq k$. With these parameters we define the {\em Pfaffian Schur process} to be the probability distribution on sequences of partitions $\lambda^1, \lambda^2, \dots, \lambda^N$, given by
\begin{equation}\label{Eq.SchurProcess}
\mathbb{P}(\lambda^1, \dots, \lambda^N) = \frac{1}{Z_N} \cdot \tau_{\lambda^1}(c) s_{\lambda^1/\lambda^2}(a_{N}) s_{\lambda^2/\lambda^3}(a_{N-1}) \cdots s_{\lambda^{N}/\emptyset}(a_1),
\end{equation}
where  $Z_N$ is a normalization constant, given explicitly in \cite[Proposition 3.2]{BR05}:
\begin{equation}\label{Eq.ZN}
   Z_N = \prod_{i = 1}^N \frac{1}{1 - ca_i} \times \prod_{1 \leq i < j \leq N} \frac{1}{1 - a_i a_j}. 
\end{equation}

Fix $J\ge 0$,  $q\in (0,1)$ and $1\le l_1<\cdots<l_J\le N$, we mostly work with the above measures when 
\begin{equation}\label{Eq.InHomogeneousParameters}
a_i=q\in(0,1) \text{ if }i\notin  \mathcal{I}:=\{l_1,\ldots,l_J\} \text{ and } a_{l_j}=\frac{1}{\zc+\sigmaq^{-1}\tilde{\alpha}_jN^{-1/3}},
\end{equation} 
in this paper. We also denote 
\begin{equation}
 \underline{\alpha}:=\min   \{\tilde{\alpha}_1,\ldots,\tilde{\alpha}_J\}.
\end{equation}
\end{definition}
\begin{remark}\label{Rem.BRSpecial} The measures in Definition \ref{Def.SchurProcess} are special cases of the Pfaffian Schur processes in \cite[Section 3]{BR05}, and correspond to setting $\rho_0^+ = c$, $\rho_i^+ = 0$ for $i = 1, \dots, N-1$, and $\rho_i^- = a_{N-i+1}$ for $i = 1, \dots, N$. 
\end{remark}
\begin{remark}
    When $J=0$, \begin{equation}\label{Eq.HomogeneousParameters}
a_1 = a_2 = \cdots = a_N = q \in (0,1),
\end{equation}
which correspond to the homogeneous case considered in Theorems \ref{Thm.Main1} and \ref{Thm.Main2}. In this case we set $\underline{\alpha}:=\infty$.
\end{remark}
\begin{lemma}\label{def: subseq}
  Suppose $(\lambda^1, \dots, \lambda^N)$ is distributed as in~(\ref{Eq.SchurProcess}). Then, for fixed indices $1 \le n_1 < n_2 < \cdots < n_k \le N$, the joint law of $\left(\lambda^{(n_1)}, \lambda^{(n_2)}, \ldots, \lambda^{(n_k)}\right)$ is given by
   \begin{equation}
\begin{aligned}
\frac{1}{Z_{N - n_1 + 1}}\,
&\tau_{\lambda^{(n_1)}}\!\left(c, a_N, \ldots, a_{N - n_1 + 2}\right)
  s_{\lambda^{(n_1)} / \lambda^{(n_2)}}\!\left(a_{N - n_1 + 1}, \ldots, a_{N - n_2 + 2}\right) \\
& \cdots 
  s_{\lambda^{(n_{k-1})} / \lambda^{(n_k)}}\!\left(a_{N - n_{k-1} + 1}, \ldots, a_{N - n_k + 2}\right)
  s_{\lambda^{(n_k)}}\!\left(a_{N - n_k + 1}, \ldots, a_1\right).
\end{aligned}
\end{equation}
\end{lemma}
\begin{proof}
    Summing over all partitions except $(\lambda^{(n_1)}, \lambda^{(n_2)}, \ldots, \lambda^{(n_k)})$ in~\eqref{Eq.SchurProcess}, 
using the skew Cauchy identity, the branching rule for Schur functions, 
and the skew Littlewood identity (see~\cite[Section~3.1.2]{BBCS18} for a summary of these identities), 
one obtains the desired result. 
\end{proof}

The following statement explains how the Pfaffian Schur processes are related to the half-space LPP model from Section \ref{Section1.1}. It is a special case of \cite[Theorem 2.7]{DY25b}. 
\begin{proposition}\label{Prop.LPPandSchur} Let $(\lambda^1, \dots,\lambda^N)$ be distributed as in (\ref{Eq.SchurProcess}) with parameters as in (\ref{Eq.InHomogeneousParameters}). Then, $(\lambda^1, \dots, \lambda^N)$ has the same distribution as $(\lambda(N,N), \lambda(N-1,N), \dots, \lambda(1,N))$ from (\ref{Eq.LPTLambdas}).
\end{proposition}
\begin{proof} From \cite[(2.10)]{DY25b} we have that almost surely $(\lambda(N,N), \lambda(N-1,N), \dots, \lambda(1,N)) = (\lambda(N,N), \lambda(N,N-1), \dots, \lambda(N,1))$. The statement now follows from \cite[Theorem 2.7]{DY25b} applied to $M = 0$ and the down-right path $\gamma = (v_0, \dots, v_N)$ with $v_i = (N, N-i)$.
\end{proof}

%
%
\subsection{Pfaffian point process structure}\label{Section2.2} In this section, we we describe how the Pfaffian Schur process can be viewed as a \emph{Pfaffian point process}. 
Throughout, we freely use the definitions and notation for Pfaffian point processes introduced in~\cite[Section~5.2]{DY25}. 
For additional background and general theory, we refer the reader to~\cite[Appendix~B]{OQR17} and~\cite{R00}.

Given two integers $a \leq b$, we  denote the set $\{a, a+1, \dots, b\}$ by $\llbracket a, b \rrbracket$. We also set $\llbracket a,b \rrbracket = \emptyset$ when $a > b$, $\llbracket a, \infty \rrbracket = \{a, a+1, a+2 , \dots \}$, $\llbracket - \infty, b\rrbracket = \{b, b-1, b-2, \dots\}$ and $\llbracket - \infty, \infty \rrbracket = \mathbb{Z}$. We fix $m \in \mathbb{N}$ and $1 \leq M_1 < M_2 < \cdots < M_m \leq N$. For a random sequence of partitions $(\lambda^1, \dots, \lambda^N)$, we define the point process $\mathfrak{S}(\lambda)$ on $ \mathbb{R}^2$ (which is supported on $\llbracket 1, m \rrbracket \times \mathbb{Z}$) by
\begin{equation}\label{Eq.PointProcessSchur}
\mathfrak{S}(\lambda)(A) = \sum_{i \geq 1} \sum_{j = 1}^m {\bf 1}\{(j, \lambda_i^{N - M_j + 1} - i) \in A\}.
\end{equation}
Below, we let $C_r$ be the positively oriented zero-centered circle of radius $r > 0$, and let $\operatorname{Mat}_2(\mathbb{C})$ be the set of $2 \times 2$ matrices with complex entries.
 \begin{proposition}\label{Prop.CorrKernel1} Assume the same notation as in Definition \ref{Def.SchurProcess} with parameters as in (\ref{Eq.InHomogeneousParameters}), and let $\mathfrak{S}(\lambda)$ be as in (\ref{Eq.PointProcessSchur}). Then, $\mathfrak{S}(\lambda)$ is a Pfaffian point process on $\mathbb{R}^2$ with reference measure given by the counting measure on $\llbracket 1, m \rrbracket \times \mathbb{Z}$ and with correlation kernel $\kgeo: (\llbracket 1, m \rrbracket \times \mathbb{Z})^2 \rightarrow \operatorname{Mat}_2(\mathbb{C})$, given as follows. For each $u,v \in \llbracket 1, m \rrbracket$ and $x,y \in \mathbb{Z}$
\begin{equation}\label{Eq.K11Geo}
\begin{split}
\kgeo_{11}(u,x;v,y) = &\frac{1}{(2\pi \im)^2} \oint_{C_{r_1}} dz \oint_{C_{r_1}} dw \frac{z-w}{(z^2 -1)(w^2 - 1)(zw - 1)} \cdot (1 - c/z) (1-c/w) \\
& \times z^{-x} w^{-y} \cdot (1- q/z)^N (1 - q/w)^N(1-qz)^{-M_u} (1- qw)^{-M_v}W(z)W(w),
\end{split}
\end{equation}
where $r_1 \in \left(1, \min\left\{ q^{-1},a_{l_1}^{-1},\cdots,a_{l_J}^{-1}\right\}\right)$ and
\begin{equation}\label{Eq.extraspike}
    W(z)=\prod_{j=1}^J \frac{(1-a_{l_j}/z)(1-qz)}{(1-q/z)(1-a_{l_j}z)}.
\end{equation}
In addition,
\begin{equation}\label{Eq.K12Geo}
\begin{split}
\kgeo_{12}(u,x;v,y) &= -\kgeo_{21}(v,y; u,x) = \frac{1}{(2\pi \im)^2} \oint_{C_{r^z_{12}}} dz \oint_{C_{r^w_{12}}} dw \frac{zw-1}{z(z-w)(z^2-1)} \cdot \frac{z-c}{w-c} \\
& \times z^{-x} w^{y} \cdot (1- q/z)^N (1 - q/w)^{-N}(1-qz)^{-M_u} (1- qw)^{M_v}\frac{W(z)}{W(w)},
\end{split}
\end{equation}
where $r_{12}^z \in \left(1, \min\left\{ q^{-1},a_{l_1}^{-1},\cdots,a_{l_J}^{-1}\right\}\right)$, $r^w_{12} > \max\left\{c,q,a_{l_1},\cdots,a_{l_J}\right\}$, and $r_{12}^w < r_{12}^z$ when $u \leq v$, while $r_{12}^z < r_{12}^w$ when $u > v$. Finally, 
\begin{equation}\label{Eq.K22Geo}
\begin{split}
\kgeo_{22}(u,x;v,y) = &\frac{1}{(2\pi \im)^2} \oint_{C_{r_2}}  dz \oint_{C_{r_2}} dw \frac{z-w}{zw -1} \cdot \frac{1}{(z-c)(w-c)}\cdot z^{x} w^{y} \\
& \times (1- q/z)^{-N} (1 - q/w)^{-N}(1-qz)^{M_u} (1- qw)^{M_v}\frac{1}{W(z)W(w)},
\end{split}
\end{equation}
where $r_2 > \max(c,q,1,a_{l_1},\cdots,a_{l_J})$. 
\end{proposition}
\begin{remark}\label{Rem.CorrKernel} Proposition \ref{Prop.CorrKernel1} is a special case of \cite[Theorem 3.3]{BR05} with the convention $i = m - u + 1$, $j = m - v + 1$, $u = x$ and $v = y$, corresponding to the choice of specializations as in Remark \ref{Rem.BRSpecial} with  $a_i$ as in \eqref{Eq.InHomogeneousParameters}. The only differences between our formulas and those in \cite[Theorem 3.3]{BR05}, see also \cite[Section 4.2]{BBCS18}, are that we changed variables $w \rightarrow 1/w$ within $K_{12}$ and $w \rightarrow 1/w, z \rightarrow 1/z$ within $K_{22}$. 
\end{remark}

%
%
\subsection{Parameter scaling for the curves}\label{Section2.3} In this section, we introduce a rescaling of the Pfaffian Schur process from Definition~\ref{Def.SchurProcess}, 
designed to capture the asymptotic behavior of 
$\{\lambda_i^j : i \ge 1,\, j \in \llbracket 1, N \rrbracket\}$.  We then derive an alternative expression for the correlation kernel from Proposition~\ref{Prop.CorrKernel1}, 
which is suitable for asymptotic analysis under this scaling.

We next introduce how we rescale our random partitions in the following definition.
\begin{definition}\label{Def.ScalingBulk} Assume the same parameters as in Definition \ref{Def.CriticalScaledLPP}. Fix $m \in \mathbb{N}$, $t_1, \dots, t_m \in \mathbb{R}$ with $t_1 < t_2 < \cdots < t_m$, and set $\mathcal{T} = \{t_1, \dots, t_m\}$. We also define for $t \in \mathcal{T}$ the quantity $T_{t} = T_t(N) =  \lfloor t N^{2/3} \rfloor$ and the lattice $\Lambda_t(N) = a_t(N) \cdot \mathbb{Z} + b_t(N)$, where 
\begin{equation}\label{Eq.LatticeBulk}
a_t(N) = (\sigmaq\zc)^{-1} N^{-1/3}, \mbox{ and } b_t(N) = (\sigmaq\zc)^{-1} N^{-1/3} \cdot \left( - \hq N - \pq T_{t} \right).
\end{equation}

We let $\mathbb{P}_N$ be the Pfaffian Schur process from Definition \ref{Def.SchurProcess} with parameters as in (\ref{Eq.InHomogeneousParameters}). Here, we assume that $N$ is sufficiently large so that 
\begin{equation}\label{Eq.LargeNBulk}
N \geq \lfloor \kappa N \rfloor +  T_{t_m}(N) >  \lfloor \kappa N \rfloor +  T_{t_{m-1}}(N) > \cdots >  \lfloor \kappa N \rfloor +  T_{t_1}(N) \geq 1.
\end{equation}
If $(\lambda^1, \dots, \lambda^N)$ have law $\mathbb{P}_N$, we define the random variables
\begin{equation}\label{Eq.XsBulk}
X_i^{j,N} = (\sigmaq\zc)^{-1} N^{-1/3} \cdot \left( \lambda_i^{N - \lfloor \kappa N \rfloor   -  T_{t_j} + 1} - \hq N - \pq T_{t_j}  - i\right) \mbox{ for } i \in \mathbb{N} \mbox{ and } j \in \llbracket 1, m \rrbracket.
\end{equation}
\end{definition}

We next introduce certain functions, contours and measures that will be used to define our alternative correlation kernel.
\begin{definition}\label{Def.SGBulk} Assume the same parameters as in Definition \ref{Def.CriticalScaledLPP}. For $z \in \mathbb{C} \setminus \{0, q, q^{-1}\}$, we introduce the functions
\begin{equation}\label{Eq.SGBulk}
    \begin{split}
        &\SFb(z) = \log(1-q/z) - \kappa \cdot \log(1-qz) - \hq \cdot \log(z), \hspace{2mm} \bar{\SFb}(z) = \SFb(z) - \SFb(\zc), \\
        &\GFb(z) = -\log(1-qz) - \pq \cdot \log(z), \hspace{2mm} \bar{\GFb}(z) = \GFb(z) - \GFb(\zc).
    \end{split}
\end{equation}
In equation (\ref{Eq.SGBulk}) and in the rest of the paper we always take the principal branch of the logarithm.
\end{definition}

\begin{definition}\label{Def.ContoursBulk} Assume the same parameters as in Definition \ref{Def.CriticalScaledLPP}. Suppose $x \in (0, q^{-1})$, $\theta \in (0, \pi)$, $R-x >|r| >0$. With this data we define the contour $C(x, \theta, R,r)$ as follows. Let $z^{\pm}$ be the points where the rays $\{x + te^{\pm \im \theta}: t \geq 0\}$ intersect the $x$-centered circle of radius $|r|$, and let $\zeta^{\pm}$ be the points where they intersect the $0$-centered circle of radius $R$. The contour $C(x, \theta, R, r)$ consists of two oriented line segments connecting $\zeta^-$ to $z^-$ and $z^+$ to $\zeta^+$, 
together with two circular arcs: the counterclockwise-oriented arc of the $0$-centered circle connecting $\zeta^+$ to $\zeta^-$, 
and the arc of the $x$-centered circle connecting $z^-$ to $z^+$, 
which is oriented clockwise if $r < 0$ and counterclockwise if $r \ge 0$.  See the left side of Figure \ref{Fig.ContoursBulk}. We also recall that $C_r$ is the positively oriented zero-centered circle of radius $r > 0$. 
\end{definition}

\begin{figure}[h]
    \centering
     \begin{tikzpicture}[scale=0.75]

        \def\tra{9} 
        \draw[->,thick,gray] (-4,0)--(4.0,0) node[right] {$\Real$};
  \draw[->,thick,gray] (0,-4)--(0,4) node[above] {$\Imag$};

  \def\x{1.2}   
  \def\R{3.4}   
  \def\r{1.0}   
  \def\th{30}   

  \node (O) at (0,0) {};
  \node (X) at (\x,0) {};
  
  \path (X) ++(\th:1) coordinate (Xp);
  \path (X) ++(-\th:1) coordinate (Xm);

  \path[name path=BigCirc]   (O) circle[radius=\R];
  \path[name path=SmallCirc] (X) circle[radius=\r];
  \path[name path=RayPlus]   (X) -- ($(X)+10*(\th:1)$);
  \path[name path=RayMinus]  (X) -- ($(X)+10*(-\th:1)$);

  \path[name intersections={of=RayPlus and BigCirc,by=ZetaPlus}];
  \path[name intersections={of=RayMinus and BigCirc,by=ZetaMinus}];
  \path[name intersections={of=RayPlus and SmallCirc,by=zPlus}];
  \path[name intersections={of=RayMinus and SmallCirc,by=zMinus}];


  \draw[dashed] (O) circle[radius=\R];
  \draw[] (zPlus) arc (30:330:\r);

  \draw[->] (ZetaMinus) -- (zMinus);
  \draw[->] (zPlus) -- (ZetaPlus);
 \draw[<->,thin] (0,0) -- ++(135:{\R})node[pos=0.5, above] {$R$};
 \draw[<->,thin] (\x,0) -- ++(90:{\r}) node[pos=0.5, left] {$|r|$};
  \draw let
      \p1 = (ZetaPlus),
      \p2 = (ZetaMinus),
      \n1 = {atan2(\y1,\x1)},                    
      \n2 = {atan2(\y2,\x2)},                    
      \n3 = {ifthenelse(\n2<\n1, \n2+360, \n2)}, 
      \n4 = {(\n1+\n3)/2}                        
    in
      [->] (ZetaPlus) arc[start angle=\n1, end angle=\n4, radius=\R];

  \draw let
      \p1 = (ZetaPlus),
      \p2 = (ZetaMinus),
      \n1 = {atan2(\y1,\x1)},
      \n2 = {atan2(\y2,\x2)},
      \n3 = {ifthenelse(\n2<\n1, \n2+360, \n2)}
    in
      (ZetaPlus) arc[start angle=\n1, end angle=\n3, radius=\R];

  
  \fill (O) circle (1.5pt) node[below left=2pt] {$0$};
  \fill (X) circle (1.5pt) node[below=2pt] {$x$};
   \fill (0.6,0) circle (1.5pt) node[below=2pt]{$z_c$};
  \fill (2.8,0) circle (1.5pt) node[below=0pt, yshift = 2pt] {$q^{-1}$};
  \fill (ZetaPlus)  circle (1.5pt) node[right=0pt] {$\zeta^{+}$};
  \fill (ZetaMinus) circle (1.5pt) node[right=0pt] {$\zeta^{-}$};
  \fill (zPlus)     circle (1.5pt) node[above = 2pt, xshift = 3pt]  {$z^{+}$};
  \fill (zMinus)    circle (1.5pt) node[below =2pt, xshift = 3pt] {$z^{-}$};

  \draw[->] (\x+1.3,0) arc (0:30:1.3);
  \node at (\x+1.5,0.4) {$\theta$};

        \draw[->, thick, gray] ({\tra -3},0)--({\tra + 4},0) node[right]{$\Real$};
        \draw[->, thick, gray] ({\tra + 0},-4)--({\tra + 0},4) node[above]{$\Imag$};
        \fill (\tra, 0) circle (1.5pt) node[below left=2pt] {$0$};

        \fill (\tra - 1.5,0) circle (1.5pt) node[above =2pt] {$a$};
        \draw[thick] ($(\tra -1.5,0)+(0.7,-0.7)$) arc[start angle=315, end angle=45, radius={sqrt(0.98)}];
        \draw[->, thick] ({\tra -0.8},0.7)--({\tra + -0.8 + 2},0.7 + 2);
        \draw[-, thick] ({\tra + -0.8 + 2},0.7 + 2)--({\tra -0.8 + 3.3},4);
         \draw[-, thick] ({\tra -0.8 + 2},-2 - 0.7)--({\tra -0.8},-0.7);
       \draw[->, thick] ({\tra + -0.8 + 3.3},-4)--({\tra -0.8 + 2},-2 - 0.7);

         \draw[<->, very thin] (\tra -1.5,0)--(\tra - 0.8, -0.7) node[pos=0.5, below left, xshift = 2pt, yshift = 2pt] {$|r|$};
        \draw[->, very thin] (\tra + 1,0) arc (0:45:2.5);
        \node at (\tra + 1.1,1) {$\varphi$};

    \end{tikzpicture} 
    \caption{The left side depicts the contours $C(x,\theta, R,r)$ with $r<0$ from Definition \ref{Def.ContoursBulk}. The right side depicts the contours $\mathcal{C}_a^{\varphi}[r]$ with $r<0$, defined above (\ref{Eq.I11Vanish}).}
    \label{Fig.ContoursBulk}
\end{figure}

\begin{definition} \label{Def:ScaledLatticeMeasures}
Fix $m \in \mathbb{N}$,  $t_1 < \cdots < t_m$, and set $\mathcal{T} = \{t_1, \dots, t_m\}$. If $\nu = (\nu_{t_1}, \dots, \nu_{t_m})$ is an $m$-tuple of locally finite measures on $\mathbb{R}$, we define the (locally finite) measure $\mu_{\mathcal{T},\nu}$ on $\mathbb{R}^2$ by
\begin{equation}\label{Eq.MuToNu}
\mu_{\mathcal{T},\nu}(A) = \sum_{t \in \mathcal{T}} \nu_t(A_{t}), \mbox{ where } A_{t} = \{ y \in \mathbb{R}: (t,y) \in A\}.
\end{equation}
\end{definition}

With the above notation in place we can state the main results of this section.

\begin{lemma}\label{Lem.PrelimitKernelsub} Assume the same notation as in Definitions \ref{Def.CriticalScaledLPP}, \ref{Def.SchurProcess}, \ref{Def.ScalingBulk} and \ref{Def.SGBulk}. In addition, we assume $c<\zc$, fix any $\theta \in (\pi/4, \pi/2)$ and $R > q^{-1}$, and let $r_1=\min\{0,2\underline{\alpha}\} $, $r_2=\min\{-1,3\underline{\alpha}\}$, we define
\begin{equation}\label{Eq.PrelimitContours}
    \begin{split}
        &\Gamma_N = C\left(\zc,\theta, R,r_1N^{-1/3} \right),\\
        & \gamma_N = C\left(\zc, 2\pi/3,\sqrt{\zc^2 + N^{-1/6} - \zc N^{-1/12}}, r_2 N^{-1/3}  \right), \\
        &\tilde{\gamma}_N = C\left(\zc, \pi/2, \sqrt{\zc^2 + N^{-1/6}}, 0  \right)
    \end{split}
\end{equation}
 as in Definition \ref{Def.ContoursBulk}. Let $M^N$ be the point process on $\mathbb{R}^2$, formed by $\{(t_j, X_i^{j,N}): i \geq 1, j \in \llbracket 1, m\rrbracket \}$. Then, for all large $N$ (depending on $q,\kappa, c, \mathcal{T}, \tilde{\alpha}_1,\ldots,\tilde{\alpha}_J $ and $\theta$) the $M^N$ is a Pfaffian point process with reference measure $\mu_{\mathcal{T},\nu(N)}$ and correlation kernel $K^N$ that are defined as follows. 

The measure $\mu_{\mathcal{T},\nu(N)}$ is as in Definition \ref{Def:ScaledLatticeMeasures} for $\nu(N) = (\nu_{t_1}(N), \dots, \nu_{t_m}(N))$, where $\nu_{t}(N)$ is $(\sigmaq\zc)^{-1} N^{-1/3}$ times the counting measure on $\Lambda_{t}(N)$. 

The correlation kernel $K^N: (\mathcal{T} \times \mathbb{R}) \times (\mathcal{T} \times \mathbb{R}) \rightarrow\operatorname{Mat}_2(\mathbb{C})$ takes the form
\begin{equation}\label{Eq:BulkKerDecomp}
\begin{split}
&K^N(s,x; t,y) = \begin{bmatrix}
    K^N_{11}(s,x;t,y) & K^N_{12}(s,x;t,y)\\
    K^N_{21}(s,x;t,y) & K^N_{22}(s,x;t,y) 
\end{bmatrix} \\
&= \begin{bmatrix}
    I^N_{11}(s,x;t,y) & I^N_{12}(s,x;t,y) + R^N_{12}(s,x;t,y) \\
    -I^N_{12}(t,y;s,x) - R^N_{12}(t,y;s,x) & I^N_{22}(s,x;t,y) 
\end{bmatrix},
\end{split}
\end{equation}
where $I^N_{ij}(s,x;t,y), R^N_{ij}(s,x;t,y)$ are defined as follows. The kernels $I^N_{ij}$ are given by
\begin{equation}\label{Eq.DefIN11Bulk}
\begin{split}
&I^N_{11}(s,x;t,y) = \frac{1}{(2\pi \im)^{2}}\oint_{\Gamma_N} dz \oint_{\Gamma_N} dw F_{11}^N(z,w) H_{11}^N(z,w)W(z)W(w) \mbox{, where }\\
& F^N_{11}(z,w) = e^{N\bar{\SFb}(z) + N\bar{\SFb}(w)} \cdot e^{T_s \bar{\GFb}(z) + T_t \bar{\GFb}(w)} \cdot e^{- \sigmaq \zc x N^{1/3} \log (z/\zc) - \sigmaq \zc y N^{1/3} \log(w/\zc)  }, \\
&H^N_{11}(z,w) = \sigmaq \zc N^{1/3} \cdot  \frac{(z-w)( 1 - c/z) (1 - c/w)(1-qz)^{\kappa N - \lfloor \kappa N \rfloor}(1-qw)^{\kappa N - \lfloor \kappa N \rfloor} }{(z^{2}-1)(w^{2}-1)(zw-1)};
 \end{split}
\end{equation}
\begin{equation}\label{Eq.DefIN12Bulk}
\begin{split}
&I^N_{12}(s,x;t,y) = \frac{1}{(2\pi \im)^{2}}\oint_{\Gamma_N} dz \oint_{\gamma_N} dw F_{12}^N(z,w) H_{12}^N(z,w) \frac{W(z)}{W(w)}\mbox{, where }\\
& F^N_{12}(z,w) = e^{N\bar{\SFb}(z) - N\bar{\SFb}(w)} \cdot e^{T_s \bar{\GFb}(z) - T_t \bar{\GFb}(w)} \cdot e^{- \sigmaq \zc x N^{1/3} \log (z/\zc) + \sigmaq \zc y N^{1/3} \log(w/\zc)  }, \\
&H^N_{12}(z,w) =  \sigmaq \zc N^{1/3} \cdot \frac{(zw - 1)(z-c)(1-qz)^{\kappa N - \lfloor \kappa N \rfloor}}{z (z-w)(z^2 - 1) (w-c)(1-qw)^{\kappa N - \lfloor \kappa N \rfloor}} ;
\end{split}
\end{equation}
\begin{equation}\label{Eq.DefIN22Bulk}
\begin{split}
&I^N_{22}(s,x;t,y) = \frac{1}{(2\pi \im)^{2}}\oint_{\gamma_N} dz \oint_{\gamma_N} dw F_{22}^N(z,w) H_{22}^N(z,w)\frac{1}{W(z)W(w)} \mbox{, where }\\
& F^N_{22}(z,w) = e^{-N\bar{\SFb}(z) - N\bar{\SFb}(w)} \cdot e^{-T_s \bar{\GFb}(z) - T_t \bar{\GFb}(w)} \cdot e^{ \sigmaq\zc x N^{1/3} \log (z/\zc) + \sigmaq \zc y N^{1/3} \log(w/\zc)  }, \\
&H^N_{22}(z,w) =   \sigmaq \zc N^{1/3}\cdot \frac{(z-w)}{(zw - 1)(z- c)(w - c)(1-qz)^{\kappa N - \lfloor \kappa N \rfloor}(1-qw)^{\kappa N - \lfloor \kappa N \rfloor}}.
\end{split}
\end{equation}
The kernel $R^N_{12}$ is given by
\begin{equation}\label{Eq.DefRN12Bulk}
\begin{split}
R^N_{12}(s,x;t,y) = &\frac{-{\bf 1}\{s > t \} \cdot \sigmaq \zc N^{1/3} }{2 \pi \im} \oint_{\tilde{\gamma}_{N}}\frac{dz}{z} e^{(T_s - T_t) \bar{\GFb}(z)} \cdot e^{ \sigmaq \zc N^{1/3}(y-x)  \log (z/\zc)}.
\end{split}
\end{equation}

\end{lemma}
\begin{proof} From Definition \ref{Def.CriticalScaledLPP} and the assumption that $c<\zc$, we have $q^{-1} > \zc>c > 1$, and so we can find $N_0$, depending on $q, \kappa, c, \underline{\alpha}$, such that for $N \geq N_0$ we have $\zc - c \geq  |r_2| N^{-1/3}$ and also for $z \in \Gamma_N$ and $w \in \gamma_N$
\begin{equation}\label{Eq.ContoursNestedBulk}
|z| \geq \zc-|r_1| N^{-1/3}> |w| \ge \zc - |r_2| N^{-1/3} > c.   
\end{equation}
Throughout the proof we assume that $N$ is sufficiently large so that $N \geq N_0$ and (\ref{Eq.LargeNBulk}) holds. Let $f: \mathbb{R} \rightarrow \mathbb{R}$ be a piece-wise linear increasing bijection such that $f(i) = t_i$ for $i \in \llbracket 1, m \rrbracket$. Define $\phi_N: \mathbb{R}^2 \rightarrow \mathbb{R}^2$ through 
$$\phi_N(s, x) = \left(f(s), (\sigmaq\zc)^{-1} N^{-1/3} \cdot \left( x- \hq N - \pq \lfloor f(s) N^{2/3} \rfloor \right) \right),$$   
and observe that $M^N = \mathfrak{S}(\lambda) \phi_N^{-1}$, where $\mathfrak{S}(\lambda)$ is as in (\ref{Eq.PointProcessSchur}) for $M_j = \lfloor \kappa N \rfloor + T_{t_j}$. It follows from Proposition \ref{Prop.CorrKernel1} and the change of variables formula \cite[Proposition 5.8(5)]{DY25} with the above $\phi_N$, \cite[Proposition 5.8(4)]{DY25} with 
\begin{equation}\label{Eq.f}
f(s,x) =  \exp \left(\sigmaq \zc x N^{1/3} \cdot \log (\zc) - T_s \GFb(\zc) - N \SFb(\zc) \right),
\end{equation}
and \cite[Proposition 5.8(4)]{DY25} with $c_1 = c_2 = (\sigmaq\zc)^{1/2} N^{1/6}$ that $M^N$ is a Pfaffian point process with reference measure $\mu_{\mathcal{T},\nu(N)}$ and correlation kernel $\tilde{K}^N: (\mathcal{T} \times \mathbb{R})^2 \rightarrow\operatorname{Mat}_2(\mathbb{C})$, given by
\begin{equation*}
\tilde{K}^N(s,x;t,y) = \begin{bmatrix} \sigmaq \zc N^{1/3} f(s,x)f(t,y) \kgeo_{11}(\tilde{s},\tilde{x}; \tilde{t},\tilde{y}) &  \sigmaq \zc N^{1/3} \frac{f(s,x)}{f(t,y)} \kgeo_{12}(\tilde{s},\tilde{x}; \tilde{t},\tilde{y}) \\ \sigmaq \zc N^{1/3} \frac{f(t,y)}{f(s,x)} \kgeo_{21}(\tilde{s},\tilde{x}; \tilde{t},\tilde{y}) & \sigmaq \zc N^{1/3} \frac{1}{f(s,x)f(t,y)} \kgeo_{22}(\tilde{s},\tilde{x}; \tilde{t},\tilde{y}) \end{bmatrix},
\end{equation*}
where $\kgeo$ is as in Proposition \ref{Prop.CorrKernel1}, $\tilde{s} = f^{-1}(s)$, $\tilde{t} = f^{-1}(t)$ and
\begin{equation*}
\tilde{x} = \hq N + \pq T_s + \sigmaq \zc N^{1/3} x, \hspace{2mm}\tilde{y} = \hq N + \pq T_t + \sigmaq \zc N^{1/3} y.
\end{equation*}
All that remains is to show that $\tilde{K}^N$ agrees with $K^N$ as in the statement of the lemma.\\

We note that the following identities hold:
\begin{equation}\label{Eq.ChangeVarsBulk}
\begin{split}
&z^{\mp \tilde{x}} (1-q/z)^{\pm N}(1-qz)^{\mp \kappa N \mp T_s} f(s,x)^{\pm 1} = e^{\pm N \bar{\SFb}(z) \pm T_s \bar{\GFb}(z) \mp \sigmaq \zc x N^{1/3} \log(z/\zc) }, \\
&w^{\mp \tilde{y}} (1-q/w)^{\pm N}(1-qw)^{\mp \kappa N \mp T_t} f(t,y)^{\pm 1} = e^{ \pm N \bar{\SFb}(w) \pm T_t \bar{\GFb}(w) \mp \sigmaq \zc y  N^{1/3} \log(w/\zc)}.
\end{split}
\end{equation}

{\bf \raggedleft Matching $K^N_{11}$ and $K^N_{22}$.} We may deform both contours $C_{r_1}$ in the definition of $\kgeo_{11}$ in Proposition \ref{Prop.CorrKernel1} to $\Gamma_N$ without crossing any of the poles of the integrand and hence without affecting the value of the integral by Cauchy's theorem.  After we perform the deformation, multiply by $\sigmaq \zc N^{1/3} f(s,x)f(t,y)$ and apply (\ref{Eq.ChangeVarsBulk})  we obtain $\tilde{K}^N_{11}( s,x; t,y) = I^N_{11}(s,x; t,y)$. For $K^N_{22}$, the deformation to $\gamma_N$ does not cross any poles of the integrand, and the resulting formula is obtained in the same manner as for $K^N_{11}$.
\\
{\bf \raggedleft Matching  $K^N_{12}$ and $K^N_{21}$.} Since $\tilde{K}^N$ and $K^N$ are both skew-symmetric, it suffices to match $K^N_{12}$. We proceed to deform $C_{r_{12}^w}$ to $\gamma_N$, and $C_{r_{12}^z}$ to $\Gamma_N$.  If $s > t$, we cross the simple pole at $z = w$. By the Residue theorem we obtain
\begin{equation}\label{eq:K12Res}
\begin{split}
&\kgeo_{12}(\tilde{s},\tilde{x}; \tilde{t},\tilde{y}) = \frac{1}{(2\pi \im)^{2}}\oint_{\Gamma_N} d  z \oint_{\gamma_N}  dw  \frac{zw-1}{z(z-w)(z^{2}-1)} \cdot \frac{z-c}{w-c} \cdot  z^{-\tilde{x} } w^{\tilde{y} }  \\
&  \times     (1-q/z)^{N} (1-q/w)^{-N}(1-qz)^{- \lfloor \kappa N \rfloor - T_s}   (1-qw)^{\lfloor \kappa N \rfloor + T_t}  \cdot \frac{W(z)}{W(w)}\\
& - \frac{{\bf 1}\{s > t\} }{2\pi \im} \oint_{\Gamma_N} \frac{dz}{z} (1-q/z)^{T_t-T_s} \cdot z^{\tilde{y}-\tilde{x}}.
\end{split}
\end{equation}
Using (\ref{Eq.ChangeVarsBulk}), one readily verifies that 
\begin{equation}\label{eq:K12Match1}
\begin{split}
&\sigmaq \zc N^{1/3}\cdot \frac{f(s,x)}{f(t,y)} \times [\mbox{lines 1 and 2 in (\ref{eq:K12Res})}] =  I_{12}^N(s,x; t,y) , \mbox{ and } \\
&\sigmaq \zc N^{1/3} \cdot \frac{f(s,x)}{f(t,y)} \times [\mbox{lines 3  in (\ref{eq:K12Res})}] = R_{12}^N(s,x;t,y).
\end{split}
\end{equation}
We mention that in the second line in (\ref{eq:K12Match1}), when matching the third line in (\ref{eq:K12Res}) with  (\ref{Eq.DefRN12Bulk}) we deformed $\Gamma_N$ to $\tilde{\gamma}_N$, which can be done without crossing any poles. Equations (\ref{eq:K12Res}) and (\ref{eq:K12Match1}) show that $\tilde{K}^N_{12}( s,x; t,y) = I^N_{12}(s,x; t,y) + R^N_{12}(s,x;t,y)$.
\end{proof}


\begin{lemma}\label{Lem.PrelimitKernelscrtical} Assume the same notation as in Definitions \ref{Def.CriticalScaledLPP}, \ref{Def.SchurProcess}, \ref{Def.ScalingBulk} and \ref{Def.SGBulk}. In addition, we assume  $c= \zc-\sigmaq^{-1}\varpi N^{-1/3}$ with $\varpi\in \mathbb{R}$ such that $-\varpi<\underline{\alpha}$, fix any $\theta \in (\pi/4, \pi/2)$ and $R > q^{-1}$, and let $r_1=-2/3\varpi+1/3\underline{\alpha}$, $r_2=-1/3\varpi+2/3\underline{\alpha}$, we define
\begin{equation}\label{Eq.PrelimitContourscritical}
    \begin{split}
        &\Gamma_N = C\left(\zc,\theta, R,r_1N^{-1/3} \right),\\
        & \gamma_N = C\left(\zc, 2\pi/3, \sqrt{\zc^2 + N^{-1/6} - \zc N^{-1/12}}, r_2 N^{-1/3}  \right), \\
        &\tilde{\gamma}_N = C\left(\zc, \pi/2, \sqrt{\zc^2 + N^{-1/6}}, 0  \right)
    \end{split}
\end{equation}
 as in Definition \ref{Def.ContoursBulk}. Let $M^N$ be the point process on $\mathbb{R}^2$, formed by $\{(t_j, X_i^{j,N}): i \geq 1, j \in \llbracket 1, m\rrbracket \}$. Then, for all large $N$ (depending on $q,\kappa, c, \mathcal{T}, \tilde{\alpha}_1,\ldots,\tilde{\alpha}_J$ and $\theta$) the $M^N$ is a Pfaffian point process with reference measure $\mu_{\mathcal{T},\nu(N)}$ and correlation kernel
\begin{equation}\label{Eq:criticalKerDecomp}
    \hat{K}^N(s,x; t,y) = \begin{bmatrix}
    N^{2/3}K^N_{11}(s,x;t,y) & K^N_{12}(s,x;t,y)\\
    K^N_{21}(s,x;t,y) & N^{-2/3}K^N_{22}(s,x;t,y) 
\end{bmatrix}
\end{equation}
where $K^N_{i,j}$ are defined as in Lemma \ref{Lem.PrelimitKernelsub} with the contours replaced by those defined in (\ref{Eq.PrelimitContourscritical}). 
\end{lemma}
\begin{proof}
The proof is nearly identical to that of Lemma~\ref{Lem.PrelimitKernelsub}, 
with the following modifications. 
Under the critical scaling $c = \zc - \sigmaq^{-1} \varpi N^{-1/3}$, 
we deform the contours as in~\eqref{Eq.PrelimitContourscritical} 
so that $\Gamma_N$ encircles $\gamma_N$, $\gamma_N$ encircles $(c, 0)$, 
and the points $\{(a_{l_j}^{-1}, 0)\}_{j=1}^J$ lie outside $\Gamma_N$. 
In addition, we multiply the function $f(s,x)$ by an extra factor of $N^{1/3}$, as in~(\ref{Eq.f}), which in turn produces the additional factor of $N^{2/3}$ in the diagonal entries.
\end{proof}

%% file: Section3.tex
%
\section{Kernel convergence}\label{Section3} The goal of this section is to establish the following statements.

\begin{proposition}\label{Prop.KernelConvsub} Assume the same notation as in Lemma \ref{Lem.PrelimitKernelsub} with $\theta = \theta_0$, $R = R_0$ as in Lemma \ref{Lem.BigContour}. If $x_N, y_N \in \mathbb{R}$ are sequences such that $\lim_{N \rightarrow \infty} x_N = x$, $\lim_{N \rightarrow \infty} y_N = y$, then for any $s,t \in \mathcal{T}$
\begin{equation}\label{Eq.KernelLimitBottom}
\begin{split}
&\lim_{N \rightarrow \infty} K_{11}^N(s,x_N;t,y_N) = 0, \hspace{2mm} \lim_{N \rightarrow \infty} K_{22}^N(s,x_N;t,y_N) =  0,\\
& \lim_{N \rightarrow \infty} \hspace{-1mm}K^N_{12}(s,x_N; t,y_N) = e^{2\fq^3s^3/3 - 2\fq^3 t^3/3 + \fq s x - \fq ty} K_{\{\tilde{\alpha}_1,\ldots,\tilde{\alpha}_J\}, \emptyset}^{\mathrm{Airy}}\left(-\fq s,  x + \fq^2 s^2 ; - \fq t, y + \fq^2 t^2 \right),
\end{split}
\end{equation}
where $K_{\mathfrak{A},\mathfrak{B}}^{\mathrm{Airy}}$ is the Airy wanderer kernel from (\ref{Eq.S1AiryKer}). 
\end{proposition}

\begin{proposition}\label{Prop.KernelConvcritical} Assume the same notation as in Lemma \ref{Lem.PrelimitKernelscrtical} with $\theta = \theta_0$, $R = R_0$ as in Lemma \ref{Lem.BigContour}. If $x_N, y_N \in \mathbb{R}$ are sequences such that $\lim_{N \rightarrow \infty} x_N = x$, $\lim_{N \rightarrow \infty} y_N = y$, then for any $s,t \in \mathcal{T}$
\begin{equation}\label{Eq.KernelLimitBottom1}
\begin{split}
&\lim_{N \rightarrow \infty} \hat{K}_{11}^N(s,x_N;t,y_N) = 0, \hspace{2mm} \lim_{N \rightarrow \infty} \hat{K}_{22}^N(s,x_N;t,y_N) =  0,\\
& \lim_{N \rightarrow \infty} \hspace{-1mm}\hat{K}^N_{12}(s,x_N; t,y_N) = e^{2\fq^3s^3/3 - 2\fq^3 t^3/3 + \fq s x - \fq ty} K_{\{\tilde{\alpha}_1,\ldots,\tilde{\alpha}_J\}, \{-\varpi\}}^{\mathrm{Airy}}\left(-\fq s,  x + \fq^2 s^2 ; - \fq t, y + \fq^2 t^2 \right),
\end{split}
\end{equation}
where $K_{\mathfrak{A},\mathfrak{B}}^{\mathrm{Airy}}$ is the Airy wanderer kernel from (\ref{Eq.S1AiryKer}). 
\end{proposition}

The proof of Proposition \ref{Prop.KernelConvsub} and \ref{Prop.KernelConvcritical} are given in Section \ref{Section3.3} and \ref{Section3.4} respectively. We state  estimates for various functions along descent contours in Section \ref{Section3.1}. In Section \ref{Section3.2} we derive suitable estimates for the functions that appear in the kernel $K^N$ in Lemma \ref{Lem.PrelimitKernelsub} and $\hat{K}^N$ in Lemma \ref{Lem.PrelimitKernelscrtical} along the contours $\Gamma_N, \gamma_N$ and $\tilde{\gamma}_N$.

%
\subsection{Preliminary estimates}\label{Section3.1} In this section we collect several analytic estimates concerning the local behavior of the functions 
$S(z)$ and $G(z)$ defined in Definition~\ref{Def.SGBulk}. 
For completeness, we restate the key results from \cite[Section~3]{DZ25} that will be used in our asymptotic analysis. 
When the proofs are short and self-contained, we include them here for the reader's convenience; 
otherwise, we refer to \cite[Section 3]{DZ25} for detailed arguments.

\begin{lemma}\label{Lem.PowerSeriesSG}  Assume the notation from Definitions \ref{Def.CriticalScaledLPP} and  \ref{Def.SGBulk}. There exist constants $\delta_0 \in (0,1)$ and $C_0 > 0$, depending on $q, \kappa$ and $c$, such that $\SFb(z), \GFb(z)$ are analytic in the disk $\{|z- \zc| < 2\delta_0\}$, and the following statements hold. If $|z - \zc| \leq \delta_0$, then
\begin{equation}\label{Eq.TaylorS1}
\left| \SFb(z) - \SFb(\zc) - (\sigmaq^3/3)(z- \zc)^3  \right| \leq C_0 |z - \zc|^4,
\end{equation}
\begin{equation}\label{Eq.TaylorG1}
\left| \GFb(z) - \GFb(\zc) - \fq\sigmaq^2 (z- \zc)^2  \right| \leq C_0 |z - \zc|^3. 
\end{equation}
\end{lemma}
\begin{proof} Note that as $q < 1 < z_c < q^{-1}$, the functions $\SFb(z), \GFb(z)$  are analytic in their respective region with $\delta_0 = (1/2) \cdot \min(1, \zc - 1, q^{-1} - \zc)$. By a direct computation we have $\SFb'(\zc) = \SFb''(\zc) = \GFb'(\zc) = 0$, and 
$$\SFb'''(\zc) = \frac{2q (q + \sqrt{\kappa})^{5}}{\kappa^{1/2} (1+q \sqrt{\kappa}) (1 - q^2)^2}= 2 \sigmaq^3, \hspace{2mm} \GFb''(\zc) = \frac{q (q + \sqrt{\kappa})^3}{\kappa (1 - q^2)^2 (1 + q \sqrt{\kappa})} = 2\fq\sigmaq^2,$$
which imply (\ref{Eq.TaylorS1}) and (\ref{Eq.TaylorG1}).
\end{proof}

\begin{lemma}\label{Lem.DecayNearCritTheta} Assume the notation in Lemma \ref{Lem.PowerSeriesSG}, and fix $\theta \in (\pi/4, \pi/2)$. There exist constants $\delta_1 \in (0,\delta_0]$ and $\epsilon_1 > 0$, depending on $q, \kappa, c, \theta$, such that for $r \in [0, \delta_1]$
\begin{equation}\label{Eq.CritDecayS1}
\Real \left [\SFb(z) - \SFb(\zc) \right] \leq -\epsilon_1 \cdot r^3, \mbox{ if } z = \zc + r e^{\pm \im \theta}, 
\end{equation}
\end{lemma}
\begin{proof} Set  $\epsilon_1 =  - \sigmaq^3/6\cos(3\theta)$. We pick $\delta_1 \in (0, \delta_0]$ small enough so that $C_0 \delta_1 \leq \epsilon_1$.
If $z = \zc + r e^{\pm \im \theta}$, we have $\Real \left[ (z- \zc)^3 \right] = r^3 \cos(3 \theta)$, so by (\ref{Eq.TaylorS1}) we conclude for $r \in [0, \delta_1]$ 
$$ \Real \left [\SFb(z) - \SFb(\zc) \right] \leq r^3 \left[(\sigmaq^3/3)\cos(3 \theta) + C_0 \delta_1\right] \leq r^3[-2\epsilon_1 + C_0\delta_1] \leq - \epsilon_1 \cdot r^3 .$$
\end{proof}

\begin{lemma}\label{Lem.DecayNearCritGen} Assume the notation in Lemma \ref{Lem.PowerSeriesSG}. There exist constants $\delta_2 \in (0,\delta_0]$ and $\epsilon_2 > 0$, depending on $q, \kappa, c$, that satisfy the following statements for $r \in [0, \delta_2]$:
\begin{equation}\label{Eq.CritGrowS1}
\Real \left [\SFb(z) - \SFb(\zc) \right] \geq \epsilon_2 \cdot r^3, \mbox{ if } z = \zc + r e^{\pm \im 2\pi/3};
\end{equation}
\begin{equation}\label{Eq.CritDecayG1}
\Real \left [\GFb(z) - \GFb(\zc) \right] \leq -\epsilon_2 \cdot r^2, \mbox{ if } z = \zc + r e^{\pm \im \pi/2}.  
\end{equation}
\end{lemma}
\begin{proof} Set $\epsilon_2 = (1/2) \min (\sigmaq^3/3, \fq \sigmaq^2 )$ and pick $\delta_2 \in (0, \delta_0]$ so that $C_0 \delta_2 \leq \epsilon_2$. If $z = \zc + r e^{\pm \im 2\pi/3}$, we have $\Real \left[ (z- \zc)^3 \right] = r^3$, and so by (\ref{Eq.TaylorS1}) we conclude for $r \in [0, \delta_2]$ 
$$ \Real \left [\SFb(z) - \SFb(\zc) \right] \geq r^3 \left[(\sigmaq^3/3) - C_0 \delta_2\right] \geq r^3 \left[2\epsilon_2 - \epsilon_2\right] \geq \epsilon_2 \cdot r^3.$$
If $z = \zc + r e^{\pm \im \pi/2}$, we have $\Real \left[ (z- \zc)^2 \right] = -r^2$, and so by (\ref{Eq.TaylorG1}) we conclude for $r \in [0, \delta_2]$ 
$$ \Real \left [\GFb(z) - \GFb(\zc) \right] \leq r^2 \left[-\fq \sigmaq^2 + C_0 \delta_2\right] \leq r^2 \left[ -2 \epsilon_2 + \epsilon_2 \right] \leq - \epsilon_2 \cdot r^2.$$
\end{proof}


\begin{lemma}\label{Lem.SmallCircleS} Assume the notation from Definitions \ref{Def.CriticalScaledLPP} and  \ref{Def.SGBulk}. If $R \in (0, \zc]$ and $z(\theta) = R e^{\pm \im \theta}$, then
\begin{equation}\label{Eq.SmallCircleS}
\frac{d}{d\theta} \Real[\SFb(z(\theta))] > 0 .
\end{equation}
\end{lemma}

\begin{proof}  Set $z(\theta) = Re^{\pm \im \theta}$, and note that 
\begin{equation}\label{Eq.CircularDerivativeS}
\frac{d}{d\theta} \Real[S(z(\theta))] = Rq \sin(\theta) \cdot\frac{R^2 q^2+1-2\cos(\theta)Rq-b^2\left( R^2 - 2\cos(\theta)Rq + q^2  \right)}{\left(R^2 - 2\cos(\theta)Rq + q^2\right)\left(R^2q^2+1-2\cos(\theta)Rq\right)}.
\end{equation}
We claim that for $R \in (0, \zc]$ we have
\begin{equation}\label{Eq.CosineUB}
\frac{R^2q^2+1-b^2\left( R^2 + q^2  \right)}{2Rq(1-b^2)} \geq 1.
\end{equation}
If (\ref{Eq.CosineUB}) holds, then we obtain for $\theta \in (0, \pi)$ that  
$$ \cos(\theta) < \frac{R^2q^2+1-b^2\left( R^2 + q^2  \right)}{2Rq(1-b^2)} \mbox{, and so }R^2 q^2+1-2\cos(\theta)Rq-b^2\left( R^2 - 2\cos(\theta)Rq + q^2  \right) > 0.$$
The last displayed equation, and (\ref{Eq.CircularDerivativeS}) then imply $\frac{d}{d\theta} \Real[S(z(\theta))] > 0$ for $\theta \in (0,\pi)$.

To see why (\ref{Eq.CosineUB}) holds, we clear denominators and see it is equivalent to 
$$(1 -qb -Rq+Rb)(1 - Rb - qR + bq) \geq 0.$$
Since $R \in (0, \zc]$, have $1 - Rb - qR + bq \geq 0$, and so it suffices to show that for $R \in [0,\zc]$
$$f(R) \geq 0, \mbox{ where }f(R) = 1 -qb -Rq+Rb.$$
The latter holds as $f(R)$ is a linear function, $f(0) = 1 -qb > 0$, and $f(\zc) = \frac{b + q - 2q^2b}{b+q} > 0$.
\end{proof}

\begin{lemma}\label{Lem.MedCircles} Assume the notation from Definitions \ref{Def.CriticalScaledLPP}  and \ref{Def.SGBulk}. If $R \in (0, \infty)$ and $z(\theta) = R e^{\pm \im \theta}$, then
\begin{equation}\label{Eq.SmallCircleG}
\frac{d}{d\theta} \Real[\GFb(z(\theta))] < 0.
\end{equation}
\end{lemma}
\begin{proof}  Set $z(\theta) = Re^{\pm \im \theta}$ and note that 
\begin{equation*}
\frac{d}{d\theta} \Real [G(z(\theta))] = -\frac{Rq \sin(\theta)}{R^2q^2 + 1 - 2Rq \cos(\theta)} < 0 \mbox{ for } \theta \in (0, \pi).
\end{equation*} 
\end{proof}
 
\begin{lemma}\cite[Lemma 3.9]{DZ25}\label{Lem.BigContour} Assume the notation from Definitions \ref{Def.CriticalScaledLPP}, \ref{Def.SGBulk} and \ref{Def.ContoursBulk}. There exist $R_0 > q^{-1}$, $\theta_0 \in (\pi/4, \pi/2)$ and a function $\psi:(0,\infty) \rightarrow (0,\infty)$, depending on $q,\kappa,c$, such that for any $\varepsilon > 0 $ 
\begin{equation}\label{Eq.DecayBigContour}
\begin{split}
\Real[\SFb(z) - \SFb(\zc)] &\leq - \psi(\varepsilon) \mbox{ if } z \in C(\zc, \theta_0, R_0,0) \mbox{ and } |z - \zc| \geq \varepsilon.
\end{split}
\end{equation}
\end{lemma}

%
\subsection{Function bounds}\label{Section3.2} 

In what follows we fix parameters as in Definition \ref{Def.CriticalScaledLPP} and $\theta_0, R_0$ as in Lemma \ref{Lem.BigContour}. In addition, we assume $c<\zc$ and work with the contours $\Gamma_N, \gamma_N, \tilde{\gamma}_N$ as in (\ref{Eq.PrelimitContours}) with $\theta = \theta_0$, $R = R_0$. We also assume that $x, y \in [-L,L]$ for a fixed $L >0$. In the inequalities below we will encounter various constants $A_i,a_i > 0$ with $A_i$ sufficiently large, and $a_i$ sufficiently small, depending on $q, \kappa, c,\mathcal{T}, \theta_0, R_0, L, \tilde{\alpha}_1,\ldots,\tilde{\alpha}_J$, we do not list this dependence explicitly. In addition, the inequalities will hold provided that $N$ is sufficiently large, depending on the same set of parameters, which we will also not mention further.  

Let $\delta_1(\theta), \epsilon_1(\theta)$ be as in Lemma \ref{Lem.DecayNearCritTheta}, $\delta_2, \epsilon_2$ be as in Lemma \ref{Lem.DecayNearCritGen} and set $\delta = \min(\delta_1(\theta_0), \delta_2), \epsilon = \min(\epsilon_1(\theta_0), \epsilon_2)$. If $z \in \Gamma_N$ and $|z - \zc| \leq \delta$, we have from Lemmas \ref{Lem.PowerSeriesSG} and \ref{Lem.DecayNearCritTheta} that
\begin{equation}\label{Eq.S1BoundZClose}
\Real[\SFb(z) - \SFb(\zc)] \leq - \epsilon |z-\zc|^3 + [\sigmaq^3/3 + \epsilon]|r_1|^3N^{-1} + C_0 |r_1|^4N^{-4/3}.
\end{equation}
If $z \in \Gamma_N$ and $|z - \zc| \geq \delta \geq  |r_1|N^{-1/3}$, we have from Lemma \ref{Lem.BigContour} that
\begin{equation}\label{Eq.S1BoundZFar}
\Real[\SFb(z) - \SFb(\zc)] \leq - \psi(\delta).
\end{equation}
If $w \in \gamma_N$ and $|w - \zc| \leq N^{-1/12} \leq \delta$, we have from Lemma \ref{Lem.DecayNearCritGen} that
\begin{equation}\label{Eq.S1BoundWClose}
\Real[\SFb(w) - \SFb(\zc)] \geq  \epsilon |w-\zc|^3. 
\end{equation}
If $w \in \gamma_N$ and $|w - \zc| \geq N^{-1/12}$, we have from Lemma \ref{Lem.SmallCircleS} that
\begin{equation}\label{Eq.S1BoundWFar}
\begin{split}
&\Real[\SFb(w) - \SFb(\zc)] \geq \Real[\SFb(\zc - e^{\pm \im 2\pi/3} N^{-1/12}) - \SFb(\zc)] \geq  \epsilon N^{-1/4}, 
\end{split}
\end{equation}
where in the last inequality we used (\ref{Eq.S1BoundWClose}). 
From Lemma \ref{Lem.PowerSeriesSG} we have for $z \in \Gamma_N \cup \gamma_N$ with $|z - \zc| \leq \delta$ 
\begin{equation}\label{Eq.G1BoundClose}
|\GFb(z) - \GFb(\zc)| \leq  \fq \sigmaq^2 |z-\zc|^2 + C_0 |z-\zc|^3.
\end{equation}
By the boundedness of the contours $\Gamma_N, \gamma_N$, we can find $A_0 > 0$, such that for $z \in \Gamma_N \cup \gamma_N$
\begin{equation}\label{Eq.G1BoundFar}
|\GFb(z) - \GFb(\zc)| \leq  A_0.
\end{equation}
By Taylor expanding the logarithm we can find $A_1 > 0$, such that for $z \in \Gamma_N \cup \gamma_N \cup \tilde{\gamma}_N$
\begin{equation}\label{Eq.G1BoundZLog}
|\log(z/\zc) | \leq  A_1 | z- \zc|.
\end{equation}

From Lemma \ref{Lem.DecayNearCritGen} we have for $z \in \tilde{\gamma}_N$ and $|z- \zc| \leq N^{-1/12} \leq \delta$
\begin{equation}\label{Eq.G1BoundZClose}
\Real[\GFb(z) - \GFb(\zc)] \leq - \epsilon |z-\zc|^2. 
\end{equation}
If $z \in \tilde{\gamma}_N$ and $|z - \zc| \geq N^{-1/12}$, we have from Lemma \ref{Lem.MedCircles} that
\begin{equation}\label{Eq.G1BoundZFar}
\begin{split}
&\Real[\GFb(z) - \GFb(\zc)] \leq \Real[\GFb(\zc + e^{\pm \im \pi/2} N^{-1/12}) - \GFb(\zc)] \leq -\epsilon N^{-1/6},
\end{split}
\end{equation}
where in the last inequality we used (\ref{Eq.G1BoundZClose}).\\

We now proceed to find suitable estimates for the functions $F^N_{ij}$, $H^N_{ij}$ and $W(z)$ from (\ref{Eq.DefIN11Bulk}), (\ref{Eq.DefIN12Bulk}), (\ref{Eq.DefIN22Bulk}) and (\ref{Eq.extraspike}). By combining (\ref{Eq.S1BoundZClose}), (\ref{Eq.S1BoundZFar}), (\ref{Eq.G1BoundFar}) and (\ref{Eq.G1BoundZLog}), we conclude that for some $A_2, a_2 > 0$ and all $z,w \in \Gamma_N$ we have 
\begin{equation}\label{Eq.F11Bound}
\begin{split}
&\left| F^N_{11}(z,w)  \right| \leq \exp \left( - \psi(\delta)N/2  \right) \mbox{ if } \max(|z-\zc|, |w-\zc|) \geq \delta, \\
&\left| F^N_{11}(z,w)  \right| \leq \exp \left( - a_2 N (|z-\zc|^3 + |w-\zc|^3) + A_2 N^{2/3} ( |z-\zc|^2 + |w-\zc|^2) + A_2  \right), \\
& \mbox{ if } \max(|z-\zc|, |w-\zc|) \leq \delta.
\end{split}
\end{equation}
By combining (\ref{Eq.S1BoundZClose}-\ref{Eq.G1BoundZLog}), we conclude that for some $A_3, a_3 > 0$ and all $z \in \Gamma_N$, $w \in \gamma_N$ we have 
\begin{equation}\label{Eq.F12Bound}
\begin{split}
&\left| F^N_{12}(z,w)  \right| \leq \exp \left( - (\epsilon/2)N^{3/4}\right) \mbox{ if } |z-\zc| \geq \delta \mbox{ or }  |w-\zc| \geq N^{-1/12}, \\
&\left| F^N_{12}(z,w)  \right| \leq \exp \left( - a_3 N (|z-\zc|^3 + |w-\zc|^3) + A_3N^{2/3} (|z-\zc|^2 + |w-\zc|^2) + A_3   \right), \\
& \mbox{ if } |z-\zc| \leq \delta \mbox{ and } |w-\zc| \leq N^{-1/12}.
\end{split}
\end{equation}
By combining (\ref{Eq.S1BoundWClose}-\ref{Eq.G1BoundZLog}), we conclude that for some $A_4, a_4 > 0$ and all $z,w \in \gamma_N$ we have 
\begin{equation}\label{Eq.F22Bound}
\begin{split}
&\left| F^N_{22}(z,w)  \right| \leq \exp \left( - (\epsilon/2)N^{3/4}\right) \mbox{ if } \max(|z-\zc|, |w-\zc|) \geq N^{-1/12}, \\
&\left| F^N_{22}(z,w)  \right| \leq \exp \left( - a_4 N (|z-\zc|^3+ |w-\zc|^3) + A_4 N^{2/3} (|z-\zc|^2 + |w-\zc|^2 ) + A_4  \right), \\
& \mbox{ if } z,w \in \max(|z-\zc|, |w-\zc|) \leq N^{-1/12}.
\end{split}
\end{equation}
As $H^N_{ij}$ are essentially rational functions and we assumed $c<\zc$, we have for some $A_5 > 0$ that
\begin{equation}\label{Eq.HijBound}
\begin{split}
&\left| H^N_{11}(z,w)  \right| \leq A_5 N^{1/3}  \mbox{ if } z, w \in \Gamma_N, \\
&\left| H^N_{12}(z,w)  \right| \leq A_5 N^{2/3}  \mbox{ if } z \in \Gamma_N \mbox{ and } w \in \gamma_N,\\
&\left| H^N_{22}(z,w)  \right| \leq A_5 N^{1/3}  \mbox{ if } z,w \in \gamma_N,
\end{split}
\end{equation}
and also 
\begin{equation}\label{Eq.HijBound2}
\begin{split}
&\left| \sigmaq \zc N^{1/3} z^{-1} \right| \leq A_5 N^{1/3}  \mbox{ if } z \in \tilde{\gamma}_N. \\
\end{split}
\end{equation}
We mention that the extra $N^{1/3}$ factor in the second line of (\ref{Eq.HijBound}) comes from the $z-w$ term in the denominator of $H^N_{12}$, for which we have $|z-w| \geq N^{-1/3}$ from the way $\Gamma_N, \gamma_N$ are defined. \\
Lastly,
\begin{equation}\label{Eq.WBound}
\begin{split}
&\left| W(z)  \right| \leq A_6 N^{J/3}  \mbox{ if } z \in \Gamma_N,\\
&\left|W(z)^{-1}  \right| \leq A_6 \mbox{ if } z \in \gamma_N,\\
& \left|\frac{W(z)}{W(w)}  \right| \leq A_6 \left|\prod_{j=1}^J\frac{w-a_{l_j}^{-1}}{z-a_{l_j}^{-1}} \right| \mbox{ if }  \mbox{ if } |z-\zc| \leq \delta \mbox{ and } |w-\zc| \leq N^{-1/12}.
\end{split}
\end{equation}

We next find suitable bounds for the functions in (\ref{Eq.DefRN12Bulk}). From (\ref{Eq.G1BoundZLog}), (\ref{Eq.G1BoundZClose}) and (\ref{Eq.G1BoundZFar}), we can find $A_7, a_7 > 0$, such that for $z \in \tilde{\gamma}_N$, $s, t \in \mathcal{T}$ with $s > t$
\begin{equation}\label{Eq.G1BoundinR}
\begin{split}
&\left| e^{(T_s - T_t) \bar{\GFb}(z) +  \sigmaq \zc N^{1/3}(y-x)  \log (z/\zc)} \right| \leq e^{- a_7 N^{1/2}} \mbox{ if } |z - \zc| \geq N^{-1/12},\\
&\left| e^{(T_s - T_t) \bar{\GFb}(z) +  \sigmaq \zc N^{1/3}(y-x)  \log (z/\zc)} \right| \leq e^{A_7 - a_7 N^{2/3}|z-\zc|^2} \mbox{ if } |z - \zc| \leq N^{-1/12}.
\end{split}
\end{equation}

%
\subsection{Proof of Proposition \ref{Prop.KernelConvsub}}\label{Section3.3} For clarity we split the proof into three steps. In the first step we show that we can truncate the contours $\gamma_N, \Gamma_N$ in the formulas for $I_{ij}^N$ in Lemma \ref{Lem.PrelimitKernelsub} without changing these functions too much, and also get some estimates for $R^N_{12}$. In the second step we prove the first line in (\ref{Eq.KernelLimitBottom}), and in the third step we prove the second line in (\ref{Eq.KernelLimitBottom}).   \\

{\bf \raggedleft Step 1.} Fix $L>0$, such that $x_N, y_N \in [-L,L]$. Let $\delta, \epsilon$ be as in the beginning of Section \ref{Section3.2}. Let $\gamma_N(0), \tilde{\gamma}_N(0)$ denote the parts of $\gamma_N, \tilde{\gamma}_N$ that are contained in the disc $\{z: |z - \zc| \leq N^{-1/12}\}$. We also denote by $\Gamma_N(0)$ the part of $\Gamma_N$ that is contained in the disc $\{z: | z -\zc| \leq \delta\}$. Let $I_{11}^{N,0}, I_{12}^{N,0}, I_{22}^{N,0},R_{12}^{N,0} $ be as in (\ref{Eq.DefIN11Bulk}), (\ref{Eq.DefIN12Bulk}), (\ref{Eq.DefIN22Bulk}), (\ref{Eq.DefRN12Bulk}) but with $\gamma_N, \Gamma_N, \tilde{\gamma}_N$ replaced with $\gamma_N(0), \Gamma_N(0), \tilde{\gamma}_N(0)$. From the first lines in (\ref{Eq.F11Bound}), (\ref{Eq.HijBound}) and (\ref{Eq.WBound})
\begin{equation}\label{Eq.TruncateI11}
\begin{split}
& \lim_{N \rightarrow \infty} \left| I^N_{11}(s,x_N; t, y_N) -  I^{N,0}_{11}(s,x_N; t, y_N) \right| \leq \lim_{N \rightarrow \infty} \|\Gamma_N\| \cdot \|\Gamma_N\| \cdot  A_5A_6^2 N^{1/3+2J/3}e^{- \psi(\delta) N/2} = 0,
\end{split}
\end{equation}
where for a contour $\gamma$, we write $\|\gamma\|$ for its arc-length. From the first line in (\ref{Eq.F12Bound}), the second line in (\ref{Eq.HijBound}) and (\ref{Eq.WBound})
\begin{equation}\label{Eq.TruncateI12}
\begin{split}
& \lim_{N \rightarrow \infty} \left| I^N_{12}(s,x_N; t, y_N) -  I^{N,0}_{12}(s,x_N; t, y_N) \right| \leq \lim_{N \rightarrow \infty} \|\Gamma_N\| \cdot \|\gamma_N\| \cdot  A_5A_6^2 N^{2/3+J/3}e^{- (\epsilon/2) N^{3/4}}= 0.
\end{split}
\end{equation}
From the first line in (\ref{Eq.F22Bound}), the third line in (\ref{Eq.HijBound}) and (\ref{Eq.WBound})
\begin{equation}\label{Eq.TruncateI22}
\begin{split}
& \lim_{N \rightarrow \infty} \left| I^N_{22}(s,x_N; t, y_N) -  I^{N,0}_{22}(s,x_N; t, y_N) \right| \leq \lim_{N \rightarrow \infty} \|\gamma_N\| \cdot \|\gamma_N\| \cdot  A_5A_6^2 N^{1/3}e^{- (\epsilon/2) N^{3/4}}= 0.
\end{split}
\end{equation}

We also note from (\ref{Eq.G1BoundZLog}), (\ref{Eq.HijBound2}) and (\ref{Eq.G1BoundinR}) 
\begin{equation}\label{Eq.R12BulkDecay}
\begin{split}
& \lim_{N \rightarrow \infty} \left| R^N_{12}(s,x_N; t, y_N) - R^{N,0}_{12} \right| \leq \lim_{N \rightarrow \infty}  \|\tilde{\gamma}_N\|  A_5N^{1/3} e^{-a_7 N^{1/2}}= 0.
\end{split}
\end{equation}

{\bf \raggedleft Step 2.} Recall from Definition \ref{Def.AiryLE} the contour $\mathcal{C}_{a}^{\varphi}=\{a+|s|e^{\mathrm{sgn}(s)\im\varphi}, s\in \mathbb{R}\}$. For $r \in \mathbb{R}$, we define the contour $\mathcal{C}_{a}^{\varphi}[r]$ so that, outside the disc $\{z : |z - a| \le |r|\}$, it agrees with $\mathcal{C}_{a}^{\varphi}$, while inside the disc it follows an arc centered at $a$ connecting the points $a + |r| e^{-\im \varphi}$ and $a + |r| e^{\im \varphi}$. When $r \ge 0$, the arc is traversed counterclockwise; when $r < 0$, it is traversed clockwise.
 see Figure \ref{Fig.ContoursBulk}. By changing variables $z = \zc + \tilde{z}N^{-1/3}$ and $w = \zc + \tilde{w}N^{-1/3}$, and applying the estimates from the second line in (\ref{Eq.F11Bound}) and the first line in (\ref{Eq.HijBound}) we obtain
\begin{equation}\label{Eq.I11Vanish}
\left|I^{N,0}_{11}(s,x_N; t, y_N) \right| \leq \frac{A_5 N^{-1/3}}{(2\pi)^2} \int_{\mathcal{C}_{0}^{\theta_0}[r_1]} |d\tilde{z}| \int_{\mathcal{C}_{0}^{\theta_0}[r_1]}|d\tilde{w}| e^{- a_2 (|\tilde{z}|^3 + |\tilde{w}|^3) + A_2 ( |\tilde{z}|^2 + |\tilde{w}|^2) + A_2},
\end{equation}
where $r_1$ is defined in Lemma \ref{Lem.PrelimitKernelsub}, and $|d\tilde{z}|, |d\tilde{w}|$ denote integration with respect to arc-length. Note that the integral on the right side of (\ref{Eq.I11Vanish}) is finite due to the cubic terms in the exponential. Consequently, the right side in (\ref{Eq.I11Vanish}) vanishes in the $N \rightarrow \infty$ limit due to the factor $N^{-1/3}$. Combining (\ref{Eq.TruncateI11}), (\ref{Eq.I11Vanish}) and the first line in (\ref{Eq.WBound}), we conclude the first limit in (\ref{Eq.KernelLimitBottom}).

Similarly to the previous paragraph, by changing variables $z = \zc +  \tilde{z}N^{-1/3}$ and $w = \zc + \tilde{w}N^{-1/3}$, and applying the estimates from the second line in (\ref{Eq.F22Bound}) and the third line in (\ref{Eq.HijBound}) we obtain
\begin{equation}\label{Eq.I22Vanish}
\left|I^{N,0}_{22}(s,x_N; t, y_N) \right| \leq \frac{A_5 N^{-1/3}}{(2\pi)^2} \int_{\mathcal{C}_{0}^{2\pi/3}[r_2]} |d\tilde{z}| \int_{\mathcal{C}_{0}^{2\pi/3}[r_2]}|d\tilde{w}| e^{- a_4 (|\tilde{z}|^3 + |\tilde{w}|^3) + A_4 ( |\tilde{z}|^2 + |\tilde{w}|^2) + A_4},
\end{equation}
which vanishes as $N \rightarrow \infty$ due to the factor $N^{-1/3}$. Combining (\ref{Eq.TruncateI22}) and (\ref{Eq.I22Vanish}), we conclude the second limit in (\ref{Eq.KernelLimitBottom}).\\

{\bf \raggedleft Step 3.} We start by computing the limit of $I^{N,0}_{12}(s,x_N; t, y_N)$. Changing variables $z = \zc + \tilde{z}N^{-1/3}$, $w = \zc + \tilde{w} N^{-1/3}$, we get
\begin{equation}\label{Eq.I12Bulk1}
\begin{split}
&I^{N,0}_{12}(s,x_N; t, y_N) =  \frac{1}{(2\pi \im)^{2}}\int_{\mathcal{C}_{0}^{\theta_0}[r_1]} d\tilde{z} \int_{\mathcal{C}_{0}^{2\pi/3}[r_2]} d\tilde{w} {\bf 1} \{|\tilde{w}| \leq N^{1/4}, |\tilde{z}| \leq \delta N^{1/3}\} \\
&\times F_{12}^N(\zc + \tilde{z}N^{-1/3},\zc + \tilde{w}N^{-1/3}) \cdot N^{-2/3} H^{N}_{12}(\zc + \tilde{z}N^{-1/3},\zc + \tilde{w}N^{-1/3})\cdot \frac{W(\zc + \tilde{z}N^{-1/3})}{W(\zc + \tilde{w}N^{-1/3})} .
\end{split}
\end{equation}
From the Taylor expansion formulas for $\SFb, \GFb$ in Lemma \ref{Lem.PowerSeriesSG} we have the pointwise limit
\begin{equation}\label{Eq.F12PointwiseBulk}
\begin{split}
\lim_{N \rightarrow \infty} F_{12}^N(\zc + \tilde{z}N^{-1/3},\zc + \tilde{w}N^{-1/3}) = e^{\sigmaq^3\tilde{z}^3/3 - \sigmaq^3\tilde{w}^3/3 +s \fq\sigmaq^2\tilde{z}^2 - t\sigmaq^2\fq \tilde{w}^2 - \sigmaq x \tilde{z} + \sigmaq y \tilde{w}}. 
\end{split}
\end{equation}
In addition, directly from the definition of $H^N_{12}$ in (\ref{Eq.DefIN12Bulk}) and $W(z)$ in (\ref{Eq.extraspike}) we get the pointwise limits
\begin{equation}\label{Eq.H12PointwiseBulk}
\begin{split}
\lim_{N \rightarrow \infty} N^{-2/3}H_{12}^N(\zc + \tilde{z}N^{-1/3},\zc + \tilde{w}N^{-1/3}) = \frac{\sigmaq}{\tilde{z} - \tilde{w}},
\end{split}
\end{equation}
\begin{equation}\label{Eq.WPointwiseBulk}
\begin{split}
\lim_{N \rightarrow \infty} \frac{W(\zc + \tilde{z}N^{-1/3})}{W(\zc + \tilde{w}N^{-1/3})} = \prod_{j=1}^J\frac{\sigmaq \tilde{w}-\tilde{\alpha}_j}{\sigmaq \tilde{z}-\tilde{\alpha}_j}.
\end{split}
\end{equation}
We also have from the second lines in (\ref{Eq.F12Bound}) and (\ref{Eq.HijBound}), and the third line in (\ref{Eq.WBound}) that the integrand in (\ref{Eq.I12Bulk1}) is bounded in absolute value by 
$$ A_5A_6\exp \left( - a_3(|\tilde{z}|^3 + |\tilde{w}|^3) + A_3 (|\tilde{z}|^2 + |\tilde{w}|^2) + A_3 \right) \left|\prod_{j=1}^J\frac{\sigmaq \tilde{w}-\tilde{\alpha}_j}{\sigmaq \tilde{z}-\tilde{\alpha}_j} \right|.$$ 
The last few observations and the dominated convergence theorem now yield
\begin{equation}\label{Eq.I12Bulk2}
\begin{split}
&\lim_{N \rightarrow \infty} I^{N,0}_{12}(s,x_N; t, y_N) = \frac{1}{(2\pi \im)^{2}}\int_{\mathcal{C}_{0}^{\theta_0}[\sigmaq r_1]} dz \int_{\mathcal{C}_{0}^{2\pi/3}[\sigmaq r_2]} dw \frac{e^{z^3/3 - w^3/3 +s \fq z^2 - t\fq w^2 - xz +  yw}}{z-w}\prod_{j=1}^J\frac{ w-\tilde{\alpha}_j}{z-\tilde{\alpha}_j},
\end{split}
\end{equation}
where we also changed variables $z = \sigmaq \tilde{z}$ and $w = \sigmaq \tilde{w}$. \\

We next compute the limit of $R^{N,0}_{12}(s,x_N; t, y_N)$. Setting $z = \zc + \tilde{z}N^{-1/3}$ gives
\begin{equation}\label{Eq.R12Bulk1}
\begin{split}
&R^{N,0}_{12}(s,x_N; t, y_N) = \frac{-{\bf 1}\{s > t \} \sigmaq \zc }{2 \pi \im} \int_{\im \mathbb{R}} \frac{d\tilde{z}}{\zc + \tilde{z}N^{-1/3}} \cdot  {\bf 1} \{|\tilde{z}| \leq N^{1/4} \} \\
&\times e^{(T_s - T_t) \bar{\GFb}(\zc + \tilde{z}N^{-1/3})} \cdot e^{ \sigmaq \zc N^{1/3}(y_N-x_N)  \log (1 + \zc^{-1} \tilde{z} N^{-1/3})}.
\end{split}
\end{equation}
From the Taylor expansion of $\GFb$ in Lemma \ref{Lem.PowerSeriesSG} the integrand in (\ref{Eq.R12Bulk1}) converges pointwise to 
$$ \frac{1}{\zc} \cdot e^{(s-t) \fq \sigmaq^2 \tilde{z}^2 + \sigmaq (y-x) \tilde{z}}.$$
We also have from  (\ref{Eq.G1BoundClose}) and (\ref{Eq.HijBound2}) that the integrand in (\ref{Eq.R12Bulk1}) is bounded in absolute value by
$$ A_5 \exp \left(A_7 - a_7|\tilde{z}|^2  \right).$$
The last few observations and the dominated convergence theorem now yield
\begin{equation}\label{Eq.R12Bulk2}
\begin{split}
&\lim_{N \rightarrow \infty} R^{N,0}_{12}(s,x_N; t, y_N) = \frac{-{\bf 1}\{s > t \}  }{2 \pi \im} \int_{\im \mathbb{R}} dz e^{(s-t) \fq z^2 + (y-x)z}
\end{split}
\end{equation}
where we also changed variables $z = \sigmaq \tilde{z}$. \\

Combining (\ref{Eq.TruncateI12}), (\ref{Eq.R12BulkDecay}), (\ref{Eq.I12Bulk2}) and (\ref{Eq.R12Bulk2}), we conclude 
\begin{equation}\label{Eq.K12LimBulk}
\begin{split}
\lim_{N \rightarrow \infty} K^N_{12} (s,x_N; t, y_N) = & -  \frac{{\bf 1}\{s > t \}  }{2 \pi \im} \int_{\im \mathbb{R}} dz e^{(s-t) \fq z^2 + (y-x)z} \\
&+ \frac{1}{(2\pi \im)^{2}}\int_{\mathcal{C}_{0}^{\theta_0}[\sigmaq r_1]} dz \int_{\mathcal{C}_{0}^{2\pi/3}[\sigmaq r_2]} dw \frac{e^{z^3/3 - w^3/3 +s \fq z^2 - t\fq w^2 - xz +  yw}}{z-w}\prod_{j=1}^J\frac{ w-\tilde{\alpha}_j}{z-\tilde{\alpha}_j}.
\end{split}
\end{equation}
What remains to prove the second line in (\ref{Eq.KernelLimitBottom}) is to verify that the right side of (\ref{Eq.K12LimBulk}) agrees with
$$A \cdot  K_{\{\tilde{\alpha}_1,\ldots,\tilde{\alpha}_J\}, \emptyset}^{\mathrm{Airy}}\left(-\fq s,  x + \fq^2 s^2 ; - \fq t, y + \fq^2 t^2 \right), \mbox{ where }  A = e^{2\fq^3s^3/3 - 2\fq^3 t^3/3 + \fq s x - \fq ty} .$$
Using the formula for $K^{\mathrm{Airy}}$ from (\ref{Eq.S1AiryKer}) with changed variables $z \rightarrow z -t_1$ and $w \rightarrow w - t_2$, we see that it suffices to show that
\begin{equation}\label{Eq.MatchBulk}
\begin{split}
& [\mbox{right side of (\ref{Eq.K12LimBulk})}] = -  A \cdot \frac{{\bf 1}\{ s > t\} }{\sqrt{4\pi \fq (s - t)}}  e^{ - \frac{(y + \fq^2 t^2 - x - \fq^2 s^2)^2}{4\fq(s - t)} - \frac{\fq(s - t)(x + y + \fq^2 s^2 + \fq^2 t^2)}{2} + \frac{\fq^3(s - t)^3}{12} }  \\
& + \frac{A}{(2\pi \im)^2} \int_{\mathcal{C}_{\alpha}^{\pi/3}} d z \int_{\mathcal{C}_{\beta}^{2\pi/3}} dw \frac{e^{(z+ \fq s)^3/3 -(x+ \fq^2 s^2)(z+ \fq s ) - (w+ \fq t)^3/3 + (y + \fq^2 t^2)(w-t_2) }}{z - w}\prod_{j=1}^J\frac{ w-\tilde{\alpha}_j}{z-\tilde{\alpha}_j}.
\end{split}
\end{equation}
where $\alpha > \beta$ and $\alpha<\underline{\alpha}$. After deforming $\mathcal{C}_{0}^{\theta_0}[\sigmaq r_1]$ and $ \mathcal{C}_{0}^{2\pi/3}[\sigmaq r_2]$ to $\mathcal{C}_{\alpha}^{\pi/3}$ and $\mathcal{C}_{\beta}^{2\pi/3}$, respectively, we see that the second line of (\ref{Eq.K12LimBulk}) agrees with the second line of (\ref{Eq.MatchBulk}). We now change variables $z = \im u$, and directly compute for $s > t$
$$\frac{1 }{2 \pi \im} \int_{\im \mathbb{R}} dz e^{(s-t) \fq z^2 + (y-x)z} =   \frac{1  }{2 \pi} \int_{ \mathbb{R}}  e^{- (s-t) \fq u^2 + \im (y-x)u} du =  \frac{1 }{\sqrt{4 \pi \fq (s-t)} }  \cdot e^{-\frac{(y-x)^2}{4\fq (s-t)}},$$
where in the last equality we used the formula for the characteristic function of a Gaussian with mean $0$ and variance $\frac{1}{2\fq (s-t)}$. Using the last identity one readily shows that the first line of (\ref{Eq.K12LimBulk}) agrees with the first line of (\ref{Eq.MatchBulk}). This completes the proof of (\ref{Eq.MatchBulk}) and hence the proposition.
%

\subsection{Proof of Proposition \ref{Prop.KernelConvcritical}}\label{Section3.4} 
In this section we assume $c= \zc-\sigmaq^{-1}\varpi N^{-1/3}$ as in Lemma \ref{Lem.PrelimitKernelscrtical} and work with the contours $\Gamma_N, \gamma_N, \tilde{\gamma}_N$ as in (\ref{Eq.PrelimitContourscritical}) with $\theta = \theta_0$, $R = R_0$ as in Lemma \ref{Lem.BigContour}. Notice that the contours $\Gamma_N, \gamma_N, \tilde{\gamma}_N$ as in (\ref{Eq.PrelimitContourscritical}) only differ from the contours $\Gamma_N, \gamma_N, \tilde{\gamma}_N$ as in (\ref{Eq.PrelimitContours}) near a small neighborhood of $(\zc,0)$ with order $N^{-1/3}$, it is not hard to see that the function bounds (\ref{Eq.S1BoundZClose})--(\ref{Eq.G1BoundinR}) derived in Section \ref{Section3.2} still hold except for the upper bounds in  (\ref{Eq.HijBound}) for $H_{i,j}^N$. Since $c= \zc-\sigmaq^{-1}\varpi N^{-1/3}$, we have  some $A_8 > 0$ that
\begin{equation}\label{Eq.HijBoundcritical}
\begin{split}
&\left| H^N_{11}(z,w)  \right| \leq A_8 N^{1/3}  \mbox{ if } z, w \in \Gamma_N, \\
&\left| H^N_{12}(z,w)  \right| \leq A_8 N  \mbox{ if } z \in \Gamma_N \mbox{ and } w \in \gamma_N,\\
&\left| H^N_{22}(z,w)  \right| \leq A_8 N  \mbox{ if } z,w \in \gamma_N.
\end{split}
\end{equation}
We mention that compared to (\ref{Eq.HijBound}), the extra $N^{1/3}$ and $N^{2/3}$ factors in the second and third line of (\ref{Eq.HijBoundcritical}) comes from the $z-c$ and $w-c$ terms in the denominators of $H_{12}^N, H_{22}^N$, for which we have $|z-c|\ge a_8 N^{1/3}$ form the way $\gamma_N$ is defined. We also note that the additional factors of $N^{2/3}$ in~(\ref{Eq:criticalKerDecomp}) are required to ensure that~(\ref{Eq.I11Vanish}) and~(\ref{Eq.I22Vanish}) still hold in this case.  \\
From this point on the proof is essentially the same as  in Section \ref{Section3.3} with the new bounds (\ref{Eq.HijBoundcritical}) and the new pointwise limit 
\begin{equation}\label{Eq.H12Pointwisecritical}
\begin{split}
\lim_{N \rightarrow \infty} N^{-2/3}H_{12}^N(\zc + \tilde{z}N^{-1/3},\zc + \tilde{w}N^{-1/3}) = \frac{\sigmaq}{\tilde{z} - \tilde{w}}\cdot\frac{\sigmaq \tilde{z}+\varpi}{\sigmaq \tilde{w}+\varpi}.
\end{split}
\end{equation}
We omit the details for the rest of proof. 

%% file: Section4.tex
%
\section{Weak convergence of the line ensembles}\label{Section4} 
We state the main results of the paper, from which Theorems~\ref{Thm.Main1} and~\ref{Thm.Main2} follow immediately as special cases when $J = 0$.
\begin{theorem}\label{Thm.generalMain}
Assume the same notation as in Definition~\ref{Def.CriticalScaledLPP}, and let $\{a_i\}$ be given in~(\ref{Eq.InHomogeneousParameters}). 
Further assume that either
\[
c \in [0, \zc) 
\quad \text{or} \quad 
c = \zc - \sigmaq^{-1}\varpi N^{-1/3}, \quad -\varpi<\underline{\alpha}.
\]
Then the sequence of random line ensembles $\mathcal{U}^N$ converges weakly,
\[
\mathcal{U}^{N} \Rightarrow \mathcal{U}^{\infty}
\quad \text{in } C\!\left(\mathbb{N} \times \mathbb{R}\right).
\]
Here $\mathcal{U}^{\infty} = \{\mathcal{U}^{\infty}_i\}_{i \ge 1}$ is the limiting line ensemble defined by
\[
\mathcal{U}_i^{\infty}(t)
= (2f)^{-1/2} \cdot
\begin{cases}
\mathcal{A}_i^{\emptyset,\{-\tilde{\alpha}_1,\ldots,-\tilde{\alpha}_J\}}\!\left(f t\right)
- f^2 t^2, 
& \text{if } c \in (0, \zc), \\[6pt]
\mathcal{A}_i^{\{\varpi\},\{-\tilde{\alpha}_1,\ldots,-\tilde{\alpha}_J\}}\!\left(f t\right)
- f^2 t^2, 
& \text{if } c = \zc - \sigmaq^{-1}\varpi N^{-1/3}.
\end{cases}
\]
In both cases, 
$\mathcal{A}^{\mathfrak{A},\mathfrak{B}}
= \{\mathcal{A}_i^{\mathfrak{A},\mathfrak{B}}\}_{i \ge 1}$ 
denotes the Airy wanderer line ensemble from Definition~\ref{Def.AiryLE}.
\end{theorem}

In Section \ref{Section4.1} we prove the finite-dimensional convergence of $(X^{j,N}_i: i\ge 1, j \in \llbracket 1, m \rrbracket)$  from Definition \ref{Def.ScalingBulk} using results from \cite{ED24a}, stated as Proposition \ref{Prop.FinitedimBulk}. Finally, in Section \ref{Section4.2} we prove Theorem \ref{Thm.generalMain}, using the previously established convergence of $(X^{j,N}_i:  i\ge 1, j \in \llbracket 1, m \rrbracket)$ and results from \cite{dimitrov2024tightness}.

%
\subsection{Finite-dimensional convergence}\label{Section4.1} As we will see later, the kernel convergence established in Proposition \ref{Prop.KernelConvsub} and \ref{Prop.KernelConvcritical} allows us to conclude that the point processes $M^N$ from Lemma \ref{Lem.PrelimitKernelsub} and \ref{Lem.PrelimitKernelscrtical} converge weakly. Our goal is then to upgrade this to finite-dimensional convergence for the vectors $(X^{j,N}_i: j \in \llbracket 1, m \rrbracket)$.
The goal of this section is to establish the following result.

\begin{proposition}\label{Prop.FinitedimBulk}
Assume the same notation as in Definition~\ref{Def.ScalingBulk} and Theorem \ref{Thm.generalMain}. 
Then the sequence of random vectors 
\[
X^N = \bigl(X_i^{j,N} : i \ge 1,\; j \in \llbracket 1, m \rrbracket \bigr)
\]
converges in the finite-dimensional sense as follows:
\[
X^N \Rightarrow
\begin{cases}
\bigl(\mathcal{A}_i^{\emptyset, \{-\tilde{\alpha}_1,\ldots,-\tilde{\alpha}_J\}}(\fq t_j)
- \fq^2 t_j^2 : i \ge 1,\; j \in \llbracket 1, m \rrbracket \bigr),
& \text{if } c < \zc, \\[8pt]
\bigl(\mathcal{A}_i^{\{\varpi\}, \{-\tilde{\alpha}_1,\ldots,-\tilde{\alpha}_J\}}(\fq t_j)
- \fq^2 t_j^2 : i \ge 1,\; j \in \llbracket 1, m \rrbracket \bigr),
& \text{if } c = \zc - \sigmaq^{-1}\varpi N^{-1/3}.
\end{cases}
\]
Here $\mathcal{A}^{\mathfrak{A},\mathfrak{B}} = \{\mathcal{A}_i^{\mathfrak{A},\mathfrak{B}}\}_{i \ge 1}$ 
denotes the Airy wanderer line ensemble from Definition~\ref{Def.AiryLE}.
\end{proposition}

\begin{proof}We provide the proof for the case \( c < \zc \); the other case follows by a similar argument and is omitted. We adopt the same notation as in Lemma \ref{Lem.PrelimitKernelsub}, where for each $s \in \mathcal{T}$, we set $\theta_s = \theta_0$, $R_s = R_0$ as in Lemma \ref{Lem.BigContour} for $\kappa = s$. For clarity, the proof is divided into three steps. In the first step, we assume that $\{X^{j,N}_i\}_{N \geq 1}$ is tight for each $i\ge 1, j \in \llbracket 1, m \rrbracket$, and conclude the proposition by applying  \cite[Proposition 2.19]{ED24a}. In the second step, we prove this tightness by assuming the sequence is tight from above -- see (\ref{Eq.TailY1}) -- and applying \cite[Proposition 2.21]{ED24a}. In the third step, we verify the tightness from above by estimating the upper-tails of the distributions of $\{X_1^{j,N}\}_{N \geq 1}$.\\

{\bf \raggedleft Step 1.} We claim that 
\begin{equation}\label{Eq.TightY}
\mbox{ the sequence } \{X_i^{j,N}\}_{N \geq 1} \mbox{ is tight for each $i\ge 1, j \in \llbracket 1, m \rrbracket$.}
\end{equation}
We prove (\ref{Eq.TightY}) in the steps below. Here, we assume its validity and complete the proof of the proposition.\\

We aim to apply  \cite[Proposition 2.19]{ED24a} with  
$$\hspace{2mm} X_i^j = \mathcal{A}^{\emptyset,\{-\tilde{\alpha}_1,\ldots,-\tilde{\alpha}_J\} }_1(\fq t_j) - \fq^2 t_j^2,  \mbox{ for }i\ge 1,  j \in \llbracket 1, m \rrbracket.$$
Note that condition (\cite[(2.25)]{ED24a}) 
\begin{equation}
 X^{j,N}_i(\omega) \geq X^{j, N}_{i+1}(\omega) \mbox{ for each $\omega \in \Omega$, $i \geq 1$, and $j \in \llbracket 1, m \rrbracket$}   
\end{equation}holds 
by the definition of $X_i^{j,N}$. In addition, the tightness of $\{ X_i^{j,N} \}_{N \geq 1}$ follows from (\ref{Eq.TightY}). Thus, it remains to show that the point processes $M^N$ from Lemma \ref{Lem.PrelimitKernelsub} converge weakly to the point process $M$ formed by $\{(t_j, \mathcal{A}^{\emptyset,\{-\tilde{\alpha}_1,\ldots,-\tilde{\alpha}_J\} }(\fq t_j) - \fq^2 t_j^2) : j \in \llbracket 1, m \rrbracket\}$.

From \cite[Lemma 5.9]{DY25} we conclude that $M$ is a Pfaffian point process on $\mathbb{R}^2$ with reference measure $\mu_{\mathcal{T}} \times \Leb$ and correlation kernel
\begin{equation}\label{Eq.LimKerFD}
K^{\mathrm{Pf}}(s,x;t,y) = \begin{bmatrix}
    0 & \frac{f(s,x)}{f(t,y)} \cdot K^{\infty}\left(s,x;t,y\right)\\
    - \frac{f(t,y)}{f(s,x)} \cdot K^{\infty}\left(t,y;s,x\right) & 0
\end{bmatrix}.
\end{equation}
Here, the functions $f(s,x)$ and $K^{\infty}\left(s,x;t,y\right)$ are given by
\begin{equation}\label{Eq.LimKerFD2}
\begin{split}
&K^{\infty}\left(s,x;t,y\right) = K_{\{\tilde{\alpha}_1,\ldots,\tilde{\alpha}_J\}, \emptyset}^{\mathrm{Airy}}\left(-\fq s,  x + \fq^2 s^2 ; - \fq t, y + \fq^2 t^2 \right), \hspace{2mm} f(s,x) = e^{2\fq^3s^3/3 + \fq s x}.
\end{split}
\end{equation}
From Lemma \ref{Lem.PrelimitKernelsub} we know (for large $N$) that $M^N$ is a Pfaffian point process on $\mathbb{R}^2$ with correlation kernel $K^N$ as in (\ref{Eq:BulkKerDecomp}) and reference measure $\mu_{\mathcal{T},\nu(N)}$. Proposition \ref{Prop.KernelConvsub} implies that for each fixed $s, t \in \mathcal{T}$, the kernels $K^{N}(s,\cdot;t,\cdot)$ converge uniformly over compact sets of $\mathbb{R}^2$ to $K^{\mathrm{Pf}}(s,\cdot;t,\cdot)$ from (\ref{Eq.LimKerFD}). Using that $K_{\{\tilde{\alpha}_1,\ldots,\tilde{\alpha}_J\}, \emptyset}^{\mathrm{Airy}}(s,x;t,y)$ is continuous in $x,y$ for fixed $s,t \in \mathcal{T}$ and hence locally bounded, it follows from \cite[Proposition 5.14]{DY25} that $M^N$ converges weakly to a Pfaffian point process $M^{\infty}$ with reference measure $\mu_{\mathcal{T}} \times \mathrm{Leb}$ and correlation kernel $K^{\mathrm{Pf}}$. \\

In view of \cite[Corollary 2.20]{ED24a}, to complete the proof of the proposition, it suffices to show that
\begin{equation}\label{Eq.MeasEqual}
M^{\infty}  \overset{d}{=}  M,
\end{equation}
where the random measure $M$ is defined by 
\begin{equation}\label{Eq.MeasEqual2}
M(A) = \sum_{i \geq 1} \sum_{j = 1}^m {\bf 1}\{(t_j, \mathcal{A}^{\emptyset,\{-\tilde{\alpha}_1,\ldots,-\tilde{\alpha}_J\} }_i(\fq t_j) - \fq^2 t_j^2) \in A\}.
\end{equation}
Using a kernel conjugation (see\cite[Proposition 5.8(4)]{DY25}), and \cite[Lemma 5.9]{DY25}, we know that $M^{\infty}$ is a determinantal point process with reference measure $\mu_{\mathcal{T}} \times \Leb$ and correlation kernel 
 \begin{equation*}
 K_{\{\tilde{\alpha}_1,\ldots,\tilde{\alpha}_J\}, \emptyset}^{\mathrm{Airy}}\left(-\fq s,  x + \fq^2 s^2 ; - \fq t, y + \fq^2 t^2 \right) = K_{\emptyset, \{-\tilde{\alpha}_1,\ldots,-\tilde{\alpha}_J\}}^{\mathrm{Airy}}\left( \fq t, y + \fq^2 t^2; \fq s,  x + \fq^2 s^2 \right),
\end{equation*}
where the equality follows directly from (\ref{Eq.S1AiryKer}) upon changing variables $z \rightarrow -w$, $w \rightarrow -z$. By kernel transposition (see \cite[Proposition 2.13(4)]{ED24a}), we conclude $M^{\infty}$ is a determinantal point process with reference measure $\mu_{\mathcal{T}} \times \Leb$ and correlation kernel
 \begin{equation}\label{Eq.LimKerFDDet}
K^{\mathrm{det}}(s,x;t,y) = K_{\emptyset, \{-\tilde{\alpha}_1,\ldots,-\tilde{\alpha}_J\}}^{\mathrm{Airy}}\left(\fq s,  x + \fq^2 s^2; \fq t,  y + \fq^2 t^2 \right).
\end{equation}

On the other hand, if $\tilde{M}$ is as in (\ref{Eq.RMS1}) with $\mathsf{S} = \{\fq t_1, \dots, \fq t_m\}$, we have that $M = \tilde{M}\phi^{-1}$, where $\phi(s,x) = (\fq^{-1} s, x - s^2)$. From \cite[Proposition 2.13(5)]{ED24a}, we conclude that $M$ is determinantal with reference measure $\mu_{\mathcal{T}} \times \Leb$ and correlation kernel $K^{\mathrm{det}}(s,x;t,y) $ as in (\ref{Eq.LimKerFDDet}). As $M$ and $M^{\infty}$ are both determinantal point processes with the same correlation kernel and reference measure, we conclude (\ref{Eq.MeasEqual}), cf. \cite[Proposition 2.13(3)]{ED24a}.\\
{\bf \raggedleft Step 2.} We claim that for each $j \in \llbracket 1, m \rrbracket$
\begin{equation}\label{Eq.TailY1}
\lim_{a \rightarrow \infty} \limsup_{N \rightarrow \infty} \mathbb{P}(X_1^{j,N} \geq a) = 0.
\end{equation}
We establish (\ref{Eq.TailY1}) in the next step. Here, we assume its validity and complete the proof of (\ref{Eq.TightY}). \\
We aim to apply \cite[Proposition 2.21]{ED24a}. Note that  condition (3) in \cite[Proposition 2.21]{ED24a} is verified by (\ref{Eq.TailY1}). It remains to check conditions (1) and (2). 

Define the point processes $M^{j,N}$ on $\mathbb{R}$ by
\begin{equation}\label{Eq.Mjn}
M^{j, N}(A) = \sum_{i \geq 1} {\bf 1}\{X_i^{j,N} \in A\}.
\end{equation}
From Lemma \ref{Lem.PrelimitKernelsub} and \cite[Lemma 5.13]{DY25} we know (for large $N$) that $M^{j,N}$ is a Pfaffian point process on $\mathbb{R}$ with reference measure $\nu_{t_j}(N)$ and correlation kernel $K^{j,N}(x,y) = K^N(t_j, x; t_j,y)$, where $K^N$ is as in (\ref{Eq:BulkKerDecomp}). From Proposition \ref{Prop.KernelConvsub} we know that $K^{j,N}(x,y)$ converges uniformly over compact sets of $\mathbb{R}^2$ to the kernel
\begin{equation}\label{Eq.KerYLimSlice}
K^{j,\infty}(x,y) = \begin{bmatrix}
    0 &  \frac{f(t_j,x)}{f(t_j,y)} \cdot K^{\infty}\left(t_j,x;t_j,y\right)\\
    - \frac{f(t_j,y)}{f(t_j,x)} \cdot K^{\infty}\left(t_j,y;t_j,x\right) & 0
\end{bmatrix},
\end{equation}
where $K^{\infty}$ is as in (\ref{Eq.LimKerFD2}). Moreover, the measures $\nu_{t_j}(N)$ converge vaguely to the Lebesgue measure on $\mathbb{R}$. Therefore, by \cite[Proposition 5.10]{DY25} we conclude that $M^{j,N}$ converge weakly to a Pfaffian point process $M^{j,\infty}$ with correlation kernel as in (\ref{Eq.KerYLimSlice}) and with reference measure $\Leb$. This verifies condition (1) in \cite[Proposition 2.21]{ED24a}. Lastly, \cite[Section 7.1]{ED24a} showed that $M^{j,\infty}$ are
infinite almost surely thus satisfies condition (2) in \cite[Proposition 2.21]{ED24a}. In conclusion, the sequence $\{X^{j,N}_i\}_{i \geq 1}$ satisfies the assumptions of \cite[Proposition 2.21]{ED24a}, from which (\ref{Eq.TightY}) follows.\\

{\bf \raggedleft Step 3.} In this step, we fix $j \in \llbracket 1, m \rrbracket$ and prove (\ref{Eq.TailY1}). If $M^{j,N}$ is as in (\ref{Eq.Mjn}), then for any $a \in \mathbb{R}$
\begin{equation}\label{Eq.TailBoundMoment}
\sum_{i \geq 1} \mathbb{P}\left(X^{j,N}_i \geq a \right) = \mathbb{E}\left[ \sum_{i \geq 1} {\bf 1}\{X_i^{j,N} \in  [a, \infty) \}   \right] = \mathbb{E}\left[M^{j,N}([a, \infty)) \right].
\end{equation}
As explained in Step 2, we have (for large $N$) that $M^{j,N}$ is a Pfaffian point process on $\mathbb{R}$ with reference measure $\nu_{t_j}(N)$ and correlation kernel $K^{j,N}(x,y) = K^N(t_j, x; t_j,y)$, where $K^N$ is as in (\ref{Eq:BulkKerDecomp}). From \cite[(2.13)]{ED24a} and \cite[(5.12)]{DY25} we know for any bounded Borel set $A \subset \mathbb{R}$
$$\mathbb{E}\left[M^{j,N}(A) \right] = \int_A K^N_{12}(t_j, x; t_j,x) \nu_{t_j}(N)(dx).$$
Setting $A = [a,b]$ and letting $b \rightarrow \infty$, we obtain by Lemma \ref{Lem.PrelimitKernelsub} and monotone convergence that
\begin{equation}\label{Eq.FactMomBulk}
\mathbb{E}\left[M^{j,N}([a, \infty)) \right] = \frac{1}{\sigmaq \zc N^{1/3}}\sum_{x \in \Lambda_{t_j}(N), x \geq a_N} K^N_{12}(t_j, x; t_j,x),
\end{equation}
where $a_N = \min\{y \in \Lambda_{t_j}(N): y \geq a\}$. Next, we notice that  
$$K^N_{12}(t_j, x; t_j,x) = I^N_{12}(t_j, x; t_j,x),$$
on the right side of (\ref{Eq.FactMomBulk}), exchange the order of the sum and the integrals in the definitions of $I^N_{12}$ from (\ref{Eq.DefIN12Bulk}), and evaluate the resulting geometric series to obtain
\begin{equation*}
\begin{split}
&\mathbb{E}\left[M^{j,N}([a, \infty)) \right] = U_N(a), \mbox{ where } U_N(a) = \frac{(\sigmaq\zc)^{-1} N^{-1/3} }{(2\pi \im)^{2}}\oint_{\Gamma_{ N}} \hspace{-3mm} dz \oint_{\gamma_{N}} dw  \frac{F_{12}^N(z,w) H_{12}^N(z,w)}{1 - w/z}\frac{W(z)}{W(w)}.
\end{split}
\end{equation*}
Here, $F^N_{12}, H^N_{12}$ are defined as in (\ref{Eq.DefIN12Bulk}) with $s =t = t_j$ and $x = y = a_N$. We mention that the two geometric series involved are absolutely convergent due to (\ref{Eq.ContoursNestedBulk}). Combining the last displayed equation with (\ref{Eq.TailBoundMoment}), we see that to prove (\ref{Eq.TailY1}), it suffices to show:
\begin{equation}\label{Eq.BulkFDRed1}
\limsup_{a \rightarrow \infty} \limsup_{N \rightarrow \infty} |U_N(a)| = 0.
\end{equation}
In the remainder of the step we verify (\ref{Eq.BulkFDRed1}) using the estimates from Section \ref{Section3.1} and \ref{Section3.2}. In the inequalities below we will encounter various constants $B_i,b_i > 0$ with $B_i$ sufficiently large, and $b_i$ sufficiently small, depending on $q, c, \mathcal{T}, \theta_0, R_0, \tilde{\alpha}_1,\ldots,\tilde{\alpha}_J$ -- we do not list this dependence explicitly. In addition, the inequalities will hold provided that $N$ is sufficiently large, depending on the same set of parameters, which we will also not mention further.\\
Let $U_N^0(a)$ be defined analogously to $U_N(a)$, but with $\Gamma_{ N},\gamma_N$ replaced with $\Gamma_{ N}(0), \gamma_N(0)$ as in Section \ref{Section3.3}. 
From (\ref{Eq.ContoursNestedBulk}), for $z \in \Gamma_{ N}$, $w \in \gamma_{N}$, and $x \geq 0$, we have
\begin{equation} \label{Eq.Ubound}
\left| \frac{1}{1 - w/z } \right| \leq R_0 N^{1/3}, \mbox{ and } \left|e^{-  \sigmaq x N^{1/3} \log (z/\zc) +  \sigmaq x N^{1/3} \log(w/\zc)  } \right|  \leq B_1 e^{-b_1 x}.
\end{equation}
Combining these bounds with (\ref{Eq.S1BoundZClose})--(\ref{Eq.G1BoundFar}), the second line in (\ref{Eq.HijBound}) and (\ref{Eq.WBound}), we obtain for some $B_2, b_2 > 0$ and all $a \geq 0$
\begin{equation*} 
\left| \frac{F_{12}^N(z,w) H_{12}^N(z,w)}{1 - w/z}\frac{W(z)}{W(w)} \right| \leq B_2 N^{2/3+J/3} e^{-b_2 N^{3/4}},  \mbox{ if } |z-\zc| \geq \delta \mbox{ or }  |w-\zc| \geq N^{-1/12},
\end{equation*}
where we recall that $F^N_{12}, H^N_{12}$ are as in (\ref{Eq.DefIN12Bulk}) with $s =t = t_j$ and $x = y = a_N \geq a$. The last bound and the bounded lengths of $\Gamma_{N}, \gamma_{N}$ imply 
\begin{equation}
    \limsup_{a \rightarrow \infty} \limsup_{N \rightarrow \infty} |U_N(a)-U_N^0(a)| = 0.
\end{equation}
Changing variables $z = \zc + \tilde{z}N^{-1/3}$, $w = \zc + \tilde{w} N^{-1/3}$, combining (\ref{Eq.S1BoundZClose}), (\ref{Eq.G1BoundClose}), the third line in (\ref{Eq.WBound}) and (\ref{Eq.Ubound}) we get 
\begin{equation}\label{Eq.U0bound}
\begin{split}
 &\left|U_N^0(a)\right|  \leq B_3 e^{-b_3 a_N}.   \int_{\mathcal{C}_{0}^{\theta_0}[r_1]} |d\tilde{z}| \int_{\mathcal{C}_{0}^{2\pi/3}[r_2]}|d\tilde{w}| \frac{e^{- a_3 (|\tilde{z}|^3 + |\tilde{w}|^3) + A_3 ( |\tilde{z}|^2 + |\tilde{w}|^2) }}{|\tilde{z}-\tilde{w}|}\left|\prod_{j=1}^J\frac{\sigmaq \tilde{w}-\tilde{\alpha}_j}{\sigmaq \tilde{z}-\tilde{\alpha}_j} \right|.
\end{split}
\end{equation}
 Note that the integral on the right side of (\ref{Eq.U0bound}) is finite due to the cubic terms in the exponential. Consequently, the right side in (\ref{Eq.U0bound}) decays exponentially in $a_N$ thus
 \begin{equation}
    \limsup_{a \rightarrow \infty} \limsup_{N \rightarrow \infty} |U_N^0(a)| = 0.
\end{equation}
This concludes the proof of (\ref{Eq.BulkFDRed1}) and hence the proposition.
\end{proof}

%
\subsection{Proof of Theorem \ref{Thm.generalMain}}\label{Section4.2} Assume the same parameters as in Definition \ref{Def.CriticalScaledLPP}, and let $\mathbb{P}_N$ be the Pfaffian Schur process from Definition \ref{Def.SchurProcess} with parameters as in (\ref{Eq.InHomogeneousParameters}). If $(\lambda^1, \dots, \lambda^N)$ is distributed according to $\mathbb{P}_N$, we define the $\mathbb{N}$-indexed geometric  line ensembles $\mathfrak{L}^N = \{L^N_i\}_{i \geq 1}$ on $\mathbb{Z}$ by
\begin{equation}\label{Eq.DLEBot}
L^{N}_i(s)  = \begin{cases} \lambda^{N - \lfloor \kappa N \rfloor - s + 1}_{i} - \lfloor  \hq N \rfloor  &\mbox{ if } i \geq 1, s + \lfloor \kappa N \rfloor  \in \llbracket 1, N \rrbracket,  \\
- \lfloor  \hq N \rfloor   &\mbox{ if } i \geq 1, s + \lfloor \kappa N \rfloor  \leq 0, \\
\lambda^1_{i} -  \lfloor \hq N \rfloor  &\mbox{ if } i \geq 1, s + \lfloor \kappa N \rfloor \geq N+1. 
\end{cases}
\end{equation}
By linear interpolation we can view $\mathfrak{L}^{N}$ as random elements in $C(\mathbb{N} \times \mathbb{R})$.  Lastly, we define the scaled ensembles $\mathcal{L}^{N} = \{\mathcal{L}^{N}_i\}_{i \geq 1} \in C(\mathbb{N} \times \mathbb{R})$ and $\mathcal{L}^{ N} = \{\mathcal{L}^{ N}_i\}_{i \geq 1} \in C(\mathbb{N} \times \mathbb{R})$ by
\begin{equation}\label{Eq.ScaledDLEBot}
\begin{split}
&\mathcal{L}^{N}_i(t) = [\pq (1+ \pq)]^{-1/2} N^{-1/3} \cdot \left(L^{N}_i( tN^{2/3}) -  \pq t N^{2/3} \right) \mbox{ for } i \geq 1, t\in \mathbb{R}.
\end{split}
\end{equation}
We mention that if $\mathcal{T} = \{t_1, \dots, t_m\}$ and $X_i^{j,N}$ are as in Definition \ref{Def.ScalingBulk}, then for all large $N$, depending on $q,\kappa,c, \tilde{\alpha}_1,\ldots,\tilde{\alpha}_J$ and $\mathcal{T}$, we have
\begin{equation}\label{Eq.LClosetoX}
\left| \mathcal{L}^{N}_i(\lfloor t_j N^{2/3} \rfloor ) - \frac{\sigmaq\zc}{\sqrt{\pq(1+\pq)}} \cdot X_{i}^{j,N} \right| \leq \frac{(i+\pq)\sigmaq\zc}{N^{1/3}\sqrt{\pq(1+\pq)}} \mbox{ for } i \geq 1, j \in \llbracket 1, m \rrbracket. 
\end{equation}

In the remainder of this section we establish Theorem \ref{Thm.generalMain}.
Combining Definition \ref{Def.CriticalScaledLPP} and Proposition \ref{Prop.LPPandSchur}, we see that $\mathcal{U}^{ N}$ has the same law as $\mathcal{L}^{ N}$ from (\ref{Eq.ScaledDLEBot}).  Consequently, it suffices to prove that 
\begin{equation}\label{Eq.ScaledTopWeakConv}
\mathcal{L}^{ N} \Rightarrow \mathcal{U}^{\infty} \in C(\mathbb{N} \times \mathbb{R}),
\end{equation}
where $\mathcal{U}^{\infty}$ is as in the statement of the theorem.\\

Fix $d > 0$, we first verify that $\mathfrak{L}^{  N}$ satisfy the hypotheses of \cite[Theorem 5.1]{dimitrov2024tightness} with $p = \pq$, $\sigma = \sigmaq$, $K = \infty$, $K_N = \infty$, $d_N = N^{2/3}$, $\hat{A}_N = \lfloor -dN^{2/3} \rfloor$, $\hat{B}_N = \lfloor  d N^{2/3} \rfloor$, $\alpha = -d$, $\beta = d$. Clearly,
$$d_N \rightarrow \infty, \hspace{2mm} \hat{A}_N/d_N \rightarrow \alpha, \hspace{2mm} \hat{B}_N/d_N \rightarrow \beta, \hspace{2mm} K_N \rightarrow K+ 1 \mbox{ as } N \rightarrow \infty,$$
verifying the first point in \cite[Theorem 5.1]{dimitrov2024tightness}. To verify the second point, we seek to show that $L^N_i(s)$
are tight for each $s \in (-d,d)$ and $i\ge 1$. However, it is implied by  Proposition \ref{Prop.FinitedimBulk} and (\ref{Eq.LClosetoX}). Lastly,  we show that ${\mathfrak{L}}^{  N}$ satisfies the interlacing Gibbs property as a $\mathbb{N}$-indexed geometric line ensemble on $\llbracket \hat{A}_N, \hat{B}_N\rrbracket$, verifying the third point in \cite[Theorem 5.1]{dimitrov2024tightness}.  We seek to verify the interlacing Gibbs property, which is equivalent to show that the  conditional law of $\left\{\mathfrak{L}^N(s): s=\hat{A}_N+1, \ldots, \hat{B}_N-1\right\}$, given that $\mathfrak{L}^N\left(\hat{A}_N\right)=\lambda^{\hat{A}_N}, \mathfrak{L}^N\left(\hat{B}_N\right)=\lambda^{\hat{B}_N}$, is uniform on sequences of all partitions subject to interlacing. This can be verified immediately with computing conditional probability using  definition \eqref{Def.SchurProcess} and Lemma \ref{def: subseq}. From \cite[Theorem 5.1]{dimitrov2024tightness}, we conclude that ${\mathcal{L}}^{  N} = \{{\mathcal{L}}^{  N}_{i}\}_{i \ge 1} \in C(\mathbb{N} \times (-d,d))$ is tight and any subsequential limit satisfies the Brownian Gibbs
property.

If $\mathcal{L}^{\infty}$ is any subsequential limit of $\mathcal{L}^{ N}$, we know from (\ref{Eq.LClosetoX}) and Proposition \ref{Prop.FinitedimBulk} that $\mathcal{L}^{\infty}$ has the same finite-dimensional distribution as $\mathcal{U}^\infty$. As finite-dimensional sets form a separating class in $C((-d,d))$, see \cite[Example 1.3]{Bill}, we conclude $\mathcal{L}^{\infty} = \mathcal{U}^\infty$. Overall, we see that the sequence $\mathcal{L}^{ N}$ is tight and each subsequential limit agrees with $\mathcal{U}^\infty$, which implies (\ref{Eq.ScaledTopWeakConv}).

%% file: PD.bib
@article{OSZ14,
author = {O'Connell, N. and Sepp{\"a}l{\"a}inen, T. and Zygouras, N.},
title = {Geometric {R}{S}{K} correspondence, {W}hittaker functions and symmetrized random polymers},
journal = {Invent. Math.},
year = 2014,
volume = {197},
issue = {2},
pages = {361--416}
}

@article{DY25,
author = {Dimitrov, E. and Yang, Z.},
title = {Half-space {A}iry line ensembles},
year = {2025},
note = {Preprint: arXiv:2505.01798},
}

@article{DY25b,
author = {Dimitrov, E. and Yang, Z.},
title = {A note on last passage percolation and {S}chur processes},
year = {2025},
note = {Preprint: arXiv:2510.04713},
}

@article{dimitrov2024tightness,
  title={Tightness for interlacing geometric random walk bridges},
  author={Dimitrov, E.},
 note={Preprint: arXiv:2410.23899},
  year={2024}
}

@article{BK08,
  author    = {Borodin, A. and Kuan, J.}, 
  title     = {Asymptotics of {P}lancherel measures for the infinite-dimensional unitary group},
  journal = {Adv. Math.},
  year      = 2008,
  volume = {219}, 
  pages = {894--931},
}

@article{SI04,
  title={Fluctuations of the one-dimensional polynuclear growth model in half-space},
  author={Sasamoto, T. and Imamura, T.},
  journal={J. Stat. Phys.},
  volume={115},
  pages={749--803},
  year={2004}
}

@article{OQR17,
author = {Ortmann, J. and Quastel, J. and Remenik, D.},
title = {A {P}faffian representation for flat {A}{S}{E}{P}},
journal = {Commun. Pure Appl. Math.},
year = 2017,
volume = {70},
number = {1},
pages = {3--89}
}

@article{Spohn,
author = {Pr{\" a}hofer, M. and Spohn, H.},
title = {Scale invariance of the {P}{N}{G} {D}roplet and the {A}iry process},
journal = {J. Stat. Phys.},
year = 2002,
volume = {108},
issue = {5},
pages = {1071--1106}
}

@article{AFM10,
author = {Adler, M. and Ferrari, P. and Van Moerbeke, P.},
title = {Airy processes with wanderers and new universality classes},
year = 2010,
journal = {Ann. Probab.},
volume = {38},
pages = {714--769},
}

@article{R00,
author = {Rains, E. M.},
title = {Correlation functions for symmetrized increasing subsequences},
year = 2000,
journal={arXiv preprint math/0006097},
}

@article{BR05,
  author    = {Borodin, A. and Rains, E.M.}, 
  title     = {Eynard-{M}ehta theorem, {S}chur process, and their {P}faffian analogs},
  journal = {J. Stat. Phys.},
  year      = 2005,
  volume = {121}, 
  pages = {291--317},
}

@article{ED24a,
author = {Dimitrov, E.},
title = {Airy wanderer line ensembles},
year = {2024},
pages = {1--76},
note = { Preprint: arXiv:2408.08445},
}

@article{DNV23,
  author    = {Dauvergne, D. and Nica, M. and Vir{\'a}g, B.}, 
  title     = {Uniform convergence to the {A}iry line ensemble},
  journal = {Ann. Inst. Henri Poincare (B) Probab. Stat.},
  year      = 2023,
  volume = {59}, 
  number = {4},
  pages = {2220--2256},
}

@article{AH23,
 AUTHOR = {Aggarwal, Amol and Huang, Jiaoyang},
     TITLE = {Edge statistics for lozenge tilings of polygons, {II}: {A}iry
              line ensemble},
   JOURNAL = {Forum Math. Pi},
  FJOURNAL = {Forum of Mathematics. Pi},
    VOLUME = {13},
      YEAR = {2025},
     PAGES = {Paper No. e2, 60},
      ISSN = {2050-5086},
   MRCLASS = {60D05 (52C23 60B20 60F05)},
  MRNUMBER = {4862833},
       DOI = {10.1017/fmp.2024.16},
       URL = {https://doi.org/10.1017/fmp.2024.16},
}

@article{TW05,
author = {Tracy, C.A. and Widom, H.},
title = {Matrix kernels for the {G}aussian orthogonal and symplectic ensembles},
journal = {Ann. Inst. Fourier},
year = 2005,
volume = {55},
number = {6},
pages = {2197--2207},
}

@article{CorHamA,
author = {Corwin, I. and Hammond, A.},
title = {Brownian {G}ibbs property for {A}iry line ensembles},
journal = {Invent. Math.},
year = 2014,
volume = {195},
pages = {441--508}
}

@book{Bill,
  author    = {Billingsley, P.}, 
  title     = {Convergence of probability measures, {S}econd Edition},
  publisher = {Academic Press, New York},
  year      = 1999,
}

@inproceedings{BBNV18,
  title={The free boundary {S}chur process and applications {I}},
  author={Betea, D. and Bouttier, J. and Nejjar, P. and Vuleti{\'c}, M.},
  booktitle={Ann. Henri Poincar{\' e}},
  volume={19},
  pages={3663--3742},
  year={2018},
  organization={Springer International Publishing}
}

@article{BBCS18,
  title={Pfaffian {S}chur processes and last passage percolation in a half-quadrant},
  author={Baik, J. and Barraquand, G. and Corwin, I. and Suidan, T.},
  journal={Ann. Probab.},
  volume={46},
  number={6},
  pages={3015--3089},
  year={2018}
}

@misc{DZ25,
      title={Curve separation in supercritical half-space last passage percolation}, 
      author={Evgeni Dimitrov and Zhengye Zhou},
      year={2025},
      eprint={2510.07508},
      archivePrefix={arXiv},
      primaryClass={math.PR},
      url={https://arxiv.org/abs/2510.07508},
 note={Preprint: arXiv:2510.07508} 
}

@article {Kurt03,
    AUTHOR = {Johansson, Kurt},
     TITLE = {Discrete polynuclear growth and determinantal processes},
   JOURNAL = {Comm. Math. Phys.},
  FJOURNAL = {Communications in Mathematical Physics},
    VOLUME = {242},
      YEAR = {2003},
    NUMBER = {1-2},
     PAGES = {277--329},
      ISSN = {0010-3616,1432-0916},
   MRCLASS = {82C22 (60F17 60K35)},
  MRNUMBER = {2018275},
MRREVIEWER = {Timo\ Sepp\"al\"ainen},
       DOI = {10.1007/s00220-003-0945-y},
       URL = {https://doi.org/10.1007/s00220-003-0945-y},
}

@article {Sod15,
     AUTHOR = {Sodin, Sasha},
     TITLE = {A limit theorem at the spectral edge for corners of
              time-dependent {W}igner matrices},
   JOURNAL = {Int. Math. Res. Not. IMRN},
  FJOURNAL = {International Mathematics Research Notices. IMRN},
      YEAR = {2015},
    NUMBER = {17},
     PAGES = {7575--7607},
      ISSN = {1073-7928,1687-0247},
   MRCLASS = {60B20},
  MRNUMBER = {3403994},
MRREVIEWER = {B\'alint\ Vet\H o},
       DOI = {10.1093/imrn/rnu180},
       URL = {https://doi.org/10.1093/imrn/rnu180},
}

@misc{IMS22,
      title={Solvable models in the KPZ class: approach through periodic and free boundary Schur measures}, 
      author={Takashi Imamura and Matteo Mucciconi and Tomohiro Sasamoto},
      year={2022},
      eprint={2204.08420},
      archivePrefix={arXiv},
      primaryClass={math.PR},
      url={https://arxiv.org/abs/2204.08420}, 
note={Preprint: arXiv:2204.08420} 
}

@article {Wu23,
    AUTHOR = {Wu, Xuan},
     TITLE = {Convergence of the {KPZ} line ensemble},
   JOURNAL = {Int. Math. Res. Not. IMRN},
  FJOURNAL = {International Mathematics Research Notices. IMRN},
      YEAR = {2023},
    NUMBER = {22},
     PAGES = {18901--18957},
      ISSN = {1073-7928,1687-0247},
   MRCLASS = {60K35 (35R60 60H15)},
  MRNUMBER = {4669793},
       DOI = {10.1093/imrn/rnac272},
       URL = {https://doi.org/10.1093/imrn/rnac272},
}

@article {DM18,
    AUTHOR = {Duse, Erik and Metcalfe, Anthony},
     TITLE = {Universal edge fluctuations of discrete interlaced particle
              systems},
   JOURNAL = {Ann. Math. Blaise Pascal},
  FJOURNAL = {Annales Math\'ematiques Blaise Pascal},
    VOLUME = {25},
      YEAR = {2018},
    NUMBER = {1},
     PAGES = {75--197},
      ISSN = {1259-1734,2118-7436},
   MRCLASS = {60B20 (82C22)},
  MRNUMBER = {3851336},
       URL = {http://ambp.cedram.org/item?id=AMBP_2018__25_1_75_0},
}

@misc{CM24,
      title={The Symplectic Schur Process}, 
      author={Cesar Cuenca and Matteo Mucciconi},
      year={2025},
      eprint={2407.02415},
      archivePrefix={arXiv},
      primaryClass={math-ph},
      url={https://arxiv.org/abs/2407.02415}, 
      note={Preprint: arXiv:2407.02415} 
}

@article {OR07,
    AUTHOR = {Okounkov, Andrei and Reshetikhin, Nicolai},
     TITLE = {Random skew plane partitions and the {P}earcey process},
   JOURNAL = {Comm. Math. Phys.},
  FJOURNAL = {Communications in Mathematical Physics},
    VOLUME = {269},
      YEAR = {2007},
    NUMBER = {3},
     PAGES = {571--609},
      ISSN = {0010-3616,1432-0916},
   MRCLASS = {60D05 (05A17 05B45 60C05 82B41)},
  MRNUMBER = {2276355},
MRREVIEWER = {Dimitri\ Petritis},
       DOI = {10.1007/s00220-006-0128-8},
       URL = {https://doi.org/10.1007/s00220-006-0128-8},
}

@article {FS03,
    AUTHOR = {Ferrari, Patrik L. and Spohn, Herbert},
     TITLE = {Step fluctuations for a faceted crystal},
   JOURNAL = {J. Statist. Phys.},
  FJOURNAL = {Journal of Statistical Physics},
    VOLUME = {113},
      YEAR = {2003},
    NUMBER = {1-2},
     PAGES = {1--46},
      ISSN = {0022-4715,1572-9613},
   MRCLASS = {82B20 (82B41 82D25)},
  MRNUMBER = {2012974},
MRREVIEWER = {Richard\ Kenyon},
       DOI = {10.1023/A:1025703819894},
       URL = {https://doi.org/10.1023/A:1025703819894},
}

@book{BMMK07,
  title     = {Discrete Orthogonal Polynomials: Asymptotics and Applications},
  series    = {Annals of Mathematics Studies},
  number    = {164},
  author    = {Baik, Jinho and Deift, Percy and McLaughlin, Kenneth T.-R. and Kriecherbauer, Thomas},
  year      = {2007},
  publisher = {Princeton University Press},
  address   = {Princeton, NJ}
}

@article{BBP05,
  AUTHOR = {Baik, Jinho and Ben Arous, G\'erard and P\'ech\'e, Sandrine},
  title   = {Phase transition of the largest eigenvalue for nonnull complex sample covariance matrices},
  journal = {Ann. Probab.},
  volume  = {33},
  number  = {5},
  pages   = {1643--1697},
  year    = {2005}
}

@article{BBCS16,
  author  = {Baik, Jinho and Barraquand, Guillaume and Corwin, Ivan and Suidan, Toufic},
  title   = {Facilitated exclusion process},
  journal = {J. Stat. Phys.},
  volume  = {165},
  number  = {6},
  pages   = {1051--1088},
  year    = {2016}
}

@article{Pet14,
  author    = {Petrov, Leonid},
  title     = {Asymptotics of random lozenge tilings via Gelfand--Tsetlin schemes},
  journal   = {Probability Theory and Related Fields},
  volume    = {160},
  number    = {3--4},
  pages     = {429--487},
  year      = {2014},
  doi       = {10.1007/s00440-013-0533-7}
}

@article {BR01a,
    AUTHOR = {Baik, Jinho and Rains, Eric M.},
     TITLE = {Algebraic aspects of increasing subsequences},
   JOURNAL = {Duke Math. J.},
  FJOURNAL = {Duke Mathematical Journal},
    VOLUME = {109},
      YEAR = {2001},
    NUMBER = {1},
     PAGES = {1--65},
      ISSN = {0012-7094,1547-7398},
   MRCLASS = {05E15 (05A15 05E05 60C05)},
  MRNUMBER = {1844203},
MRREVIEWER = {Ira\ Gessel},
       DOI = {10.1215/S0012-7094-01-10911-3},
       URL = {https://doi.org/10.1215/S0012-7094-01-10911-3},
}

@article {BR01b,
    AUTHOR = {Baik, Jinho and Rains, Eric M.},
     TITLE = {The asymptotics of monotone subsequences of involutions},
   JOURNAL = {Duke Math. J.},
  FJOURNAL = {Duke Mathematical Journal},
    VOLUME = {109},
      YEAR = {2001},
    NUMBER = {2},
     PAGES = {205--281},
      ISSN = {0012-7094,1547-7398},
   MRCLASS = {60C05 (05E10 45E10)},
  MRNUMBER = {1845180},
MRREVIEWER = {Bernhard\ Gittenberger},
       DOI = {10.1215/S0012-7094-01-10921-6},
       URL = {https://doi.org/10.1215/S0012-7094-01-10921-6},
}

@incollection {BR01c,
    AUTHOR = {Baik, Jinho and Rains, Eric M.},
     TITLE = {Symmetrized random permutations},
 BOOKTITLE = {Random matrix models and their applications},
    SERIES = {Math. Sci. Res. Inst. Publ.},
    VOLUME = {40},
     PAGES = {1--19},
 PUBLISHER = {Cambridge Univ. Press, Cambridge},
      YEAR = {2001},
      ISBN = {0-521-80209-1},
   MRCLASS = {82B41 (60C05 82B05 82B43)},
  MRNUMBER = {1842780},
MRREVIEWER = {Estelle\ L.\ Basor},
       DOI = {10.2977/prims/1145475964},
       URL = {https://doi.org/10.2977/prims/1145475964},
}

@article {BP08,
    AUTHOR = {Borodin, Alexei and P\'ech\'e, Sandrine},
     TITLE = {Airy kernel with two sets of parameters in directed
              percolation and random matrix theory},
   JOURNAL = {J. Stat. Phys.},
  FJOURNAL = {Journal of Statistical Physics},
    VOLUME = {132},
      YEAR = {2008},
    NUMBER = {2},
     PAGES = {275--290},
      ISSN = {0022-4715,1572-9613},
   MRCLASS = {82B43 (60K35 82B44)},
  MRNUMBER = {2415103},
MRREVIEWER = {Akira\ Sakai},
       DOI = {10.1007/s10955-008-9553-8},
       URL = {https://doi.org/10.1007/s10955-008-9553-8},
}

@article {BW23,
    AUTHOR = {Barraquand, Guillaume and Wang, Shouda},
     TITLE = {An identity in distribution between full-space and half-space
              log-gamma polymers},
   JOURNAL = {Int. Math. Res. Not. IMRN},
  FJOURNAL = {International Mathematics Research Notices. IMRN},
      YEAR = {2023},
    NUMBER = {14},
     PAGES = {11877--11929},
      ISSN = {1073-7928,1687-0247},
   MRCLASS = {60K35 (82D60)},
  MRNUMBER = {4615220},
MRREVIEWER = {Simone\ Baldassarri},
       DOI = {10.1093/imrn/rnac132},
       URL = {https://doi.org/10.1093/imrn/rnac132},
}
